\documentclass[10pt]{article}
\usepackage[francais,english]{babel}
\usepackage{amsmath}
\usepackage{amsfonts}
\usepackage{amssymb}
\usepackage{amsthm}
\usepackage{graphics}
\usepackage{amscd}
\usepackage{epsfig}
\usepackage{color}
\newcommand{\R}{\mathbb{R}}
\newcommand{\N}{\mathbb{N}}

\newcommand{\C}{\mathbb{C}}

\newcommand{\E}{\mathbb{E}}

\newcommand{\mc}{\mathcal}

\newcommand{\eps}{\varepsilon}
\newcommand{\ind}{{\bf 1}}
\renewcommand{\P}{\mathbb{P}}

\renewcommand{\H}{\mathbb{H}}

\DeclareMathOperator{\diam}{diam}

\DeclareMathOperator{\SLE}{SLE}

\DeclareMathOperator{\capp}{cap}
\DeclareMathOperator{\Beta}{Beta}
\DeclareMathOperator{\BES}{BES}

\title{Excursion decompositions for $\SLE$ and Watts' crossing formula}
\author{Julien Dub\'edat\footnote{Universit\'e Paris-Sud}}
\newtheorem{thm}{Theorem}

\newtheorem{Prop}[thm]{Proposition}
\newtheorem{Lem}[thm]{Lemma}
\newtheorem{Cor}[thm]{Corollary}

\begin{document}
\maketitle
\begin{abstract}
It is known that Schramm-Loewner Evolutions (SLEs) have a.s. frontier
points if $\kappa>4$ and a.s. cutpoints if $4<\kappa<8$. If
$\kappa>4$, an appropriate version of $\SLE(\kappa)$ has a renewal
property: it starts afresh after visiting its frontier. Thus one can
give an excursion decomposition for this particular $\SLE(\kappa)$
``away from its frontier''. For $4<\kappa<8$, there is a two-sided analogue
of this situation: a particular version of $\SLE(\kappa)$ has a
renewal property w.r.t its cutpoints; one studies excursion decompositions of this
$\SLE$ ``away from its cutpoints''. For $\kappa=6$, this overlaps
Vir\'ag's results on ``Brownian beads''. As a by-product of this
construction, one proves Watts' formula, which describes the
probability of a double crossing in a rectangle for critical plane percolation.
\end{abstract}

Schramm-Loewner Evolutions (SLEs) are a family of growth processes in
simply connected plane domains. Their conformal invariance properties
make them natural candidates to describe the scaling limit of critical
plane systems, that are generally conjectured to converge to a
conformally invariant limit. This convergence has been rigorously
established in several cases: for instance, the scaling limit of the
Loop-Erased Random Walk (resp. the Peano curve of the Uniform Spanning Tree)
is $\SLE(2)$ (resp. $\SLE(8)$\,) (see \cite{LSW2}), and the scaling
limit of critical percolation interfaces is $\SLE(6)$ (see
\cite{Sm1}).

The qualitative features of $\SLE$ depend crucially on the value of
the $\kappa$ parameter. The growth process is generated by a
continuous path, the trace (see \cite{RS01}). This path is a.s. simple
if $\kappa\leq 4$; if $\kappa>4$, it is no longer the case, and $\SLE$
has a non-trivial frontier (\cite{RS01}). Furthermore, $\SLE$ has cutpoints if
$4<\kappa<8$ (see \cite{Bef}). 

In \cite{Vir}, Vir\'ag shows that the Brownian Excursion (Brownian
motion in the upper
half-plane, started from 0 and conditioned not to hit the real line
again) can be decomposed in ``beads'', i.e. portions of the Brownian excursion between
two successive cutpoints. This decomposition can be phrased in terms
similar to It\^o's theory of Brownian excursions.

For $\SLE$, one also has a Markov property and conformal invariance,
so it is quite natural to look for similar decompositions w.r.t. loci
with an intrinsic geometrical definition. We will see that such
decompositions exist (for suitably conditioned $\SLE$s) for frontier
points and cutpoints. 

While considering restriction formulas in \cite{D4}, in relation with duality
conjectures, it appeared that a particular version of $\SLE(\kappa)$,
namely $\SLE(\kappa,\kappa-4)$, played a special role. For
$\kappa>4$, $\SLE(\kappa)$ a.s. swallows any real point;
and $\SLE_{(0,0^+)}(\kappa,\kappa-4)$ can be viewed as an $\SLE(\kappa)$
``conditioned'' not to hit the positive half-line. As this event has
zero probability, a little care (and precision) is required. This
process, run until infinity, has a right-boundary that is a simple
path connecting $0$ and $\infty$ in $\H$; this right-boundary is
conjectured to be identical in law to an $\SLE(\kappa',\kappa'/2-2)$,
where $\kappa\kappa'=16$. 

Looking at the right-boundary as a path, it is quite natural to
consider the conditional law of $\SLE(\kappa,\kappa-4)$ given an
initial portion of its right-boundary. A point on this (final)
right-boundary will never be swallowed again; so the conditional
$\SLE(\kappa,\kappa-4)$ avoids the already completed right-boundary,
and also avoids the positive half-line. It turns out that, after
taking the image under a conformal equivalence, the completed
right-boundary and the positive half-line play exactly the same
role. This suggests that the future of $\SLE(\kappa,\kappa-4)$ after a frontier time
(i.e. a time at which the trace is on the {\em final} right-boundary)
is again $\SLE(\kappa,\kappa-4)$ in the remaining domain, a renewal
property similar to the Markov property of $\SLE$. We will prove that one
can define a ``local time'' for time spent by $\SLE(\kappa,\kappa-4)$
on its right-boundary; subordinating the $\SLE(\kappa,\kappa-4)$ by the
inverse of this local time (which is a stable subordinator with index
$(1/2+2/\kappa)$\,), one gets a hull-valued L\'evy process. We will also
give an excursion decomposition with respect to the right-boundary
(i.e. we will describe the law of $\SLE(\kappa,\kappa-4)$ between 
two successive frontier times). As a by-product of this construction,
a one-parameter extension of Pitman's $(2M-X)$ theorem is derived.

As frontier points, cutpoints form a locus with an intrinsic geometrical
definition. As pointed out by Vir\'ag, the Brownian excursion
(i.e. planar Brownian motion started from 0 and conditioned not to
visit the lower half-plane) can be decomposed into ``beads'' (portions
of the Brownian excursion between successive cutpoints). This
decomposition relies on the conformal invariance of Brownian
motion, and the Strong Markov property. 

When $4<\kappa<8$, $\SLE$ is known to have cutpoints (the corresponding
cut-times have a.s. Hausdorff dimension $(2-\kappa/4)$, see
\cite{Bef}). By analogy with the one-sided construction, where
$\SLE(\kappa,\kappa-4)$ turned out to be well suited to the study of the
right-boundary, the structure of cut-times for
$\SLE_{(0,0^-,0^+)}(\kappa,\kappa-4,\kappa-4)$ is particularly
nice. This is due to the fact that
$\SLE_{(0,0^-,0^+)}(\kappa,\kappa-4,\kappa-4)$ can be viewed as
$\SLE(\kappa)$ ``conditioned not to hit $\R^*$'', via an appropriate
conditioning procedure. The construction of the local time of
cutpoints is
more involved than that of ``frontier local time'' (based on local
times for Bessel processes), but it can be carried out
explicitly. Subordinating the $\SLE(\kappa,\kappa-4,\kappa-4)$ process
by the inverse of this local time leads also to a
hull-valued L\'evy process; an It\^o-type result for the ``bead
process'' is also derived. In the particular case $\kappa=6$, as the Brownian excursion
and $\SLE(6,2,2)$ can be seen as realizations of the (unique)
restriction measure with exponent 1, these results partially overlap
those in \cite{Vir}.
  
In the first section, we derive some properties of $\SLE(\kappa,\kappa-4)$ and
  $\SLE(\kappa,\kappa-4,\kappa-4)$ processes that will be used later,
  using mainly stochastic calculus. The second section is devoted to
the study of $\SLE(\kappa,\kappa-4)$ in relation with its
right-boundary. We define a Markov process measuring the distance to the
boundary; the local time at 0 of this process defines
a ``frontier local time''.  The third section develops
analogous results in the two-sided case (i.e. in relation with
cutpoints), beginning with a systematic reinterpretation of the one-sided
situation in terms of Doob $h$-transforms/Girsanov densities. Though
the line of reasoning is essentially parallel, details are much more
intricate in the two-sided case. In the last section, we use some of
the previous results on $\SLE(6,2,2)$ to prove Watts' formula, which
describes the probability of the existence of a double crossing in a
rectangle, in the scaling limit of critical percolation.

{\bf Acknowledgments.} I wish to thank Wendelin Werner for his help
and advice along the preparation of this paper, as well as Marc Yor
for stimulating conversations.

%We have previously considered the capacity of the right boundary of a
%hull; a clear shortcoming of this notion is that it is not additive
%under concatenation. In this section we will define in the particular
%case of $\SLE(\kappa,\kappa-4)$ a measure of the ``length'' of the hull
%right boundary which is additive for concatenation; actually this ``length''
%depends on the whole Loewner path, not just on the boundary.

\section{Introduction and notations}

In this section, we briefly recall the definition and elementary
properties of $\SLE$ processes, and collect first properties of
$\SLE(\kappa,\kappa-4)$ (resp. $\SLE(\kappa,\kappa-4,\kappa-4)$) processes.

\subsection{Chordal $\SLE$ and $\SLE(\kappa,\underline\rho)$ processes}

For general background on $\SLE$ processes, introduced by Oded Schramm
in \cite{S0}, see \cite{RS01,W1}. In this article, we will consider only
chordal $\SLE$s, i.e. $\SLE$s that grow from one boundary point to
another boundary point in a simply connected plane domain.

Chordal Loewner equations are a device that encode a growth process in a
plane simply connected domain by a real-valued process. More precisely,
consider the upper half-plane $\H=\{z:\Im z>0\}$ (without loss of
generality since, by Riemann's mapping
theorem, any simply connected domain other than $\C$ is conformally
equivalent to $\H$), a function $w:\R^+\mapsto\R$, and the family of
ODEs:
$$\partial_t g_t(z)=\frac 2{g_t(z)-w_t}$$
with initial condition $g_0(z)=z$, $z\in\H$. The solution of any
of these ODEs is defined up to explosion time $\tau_z$ (possibly
infinite). Then, for $t\geq 0$, let
$$K_t=\overline{\{z\in\H :\tau_z<t\}}.$$
The increasing family $(K_t)$ of compact subsets of $\overline\H$ is
such that 
$$g_t(z)=z+\frac{2t}z+O(z^{-2})$$
is the unique conformal equivalence $\H\setminus K_t\rightarrow\H$
with asymptotic expansion at infinity $g_t(z)=z+o(1)$ (hydrodynamic
normalization). The coefficient $2t$ in this expansion is by definition
the half-plane capacity $\capp(K_t)$ of $K_t$. This construction sets up a
bijection between real-valued continuous processes and increasing
families of compact sets $(K_t)$ such that $\capp(K_t)=2t$, under a
``local growth'' condition (see e.g. \cite{W1}).

If $W/\sqrt\kappa$ is a standard (real) Brownian motion, the
associated (random) families of hulls $(K_t)$ and conformal equivalences
$(g_t)$ define the chordal Schramm-Loewner Evolution in $(\H,0,\infty)$
with parameter $\kappa$, in short $\SLE(\kappa)$. Note that, as a
consequence of Brownian scaling, $(\lambda^{-1}K_{\lambda^2t})_{t\geq
  0}$ as the same law as $(K_t)$ for any positive $\lambda$. Since the
dilatations $z\mapsto\lambda z$, $\lambda>0$, are the only conformal
automorphisms of $(\H,0,\infty)$, one can define $\SLE(\kappa)$ in any
simply connected domain $(D,a,b)$, where $a$ and $b$ are two distinct
boundary points, as the image of chordal $\SLE(\kappa)$ in
$(\H,0,\infty)$ as defined above under any conformal equivalence
between $(\H,0,\infty)$ and $(D,a,b)$. 

If $(K_t)$ is a Loewner chain associated with an $\SLE(\kappa)$
process, then for any $s\geq 0$, $(g_s(K_{t+s}\setminus K_s)-W_s)_t$
defines a chordal $\SLE(\kappa)$ process in $(\H,0,\infty)$
independent from $(K_t)_{t\leq s}$. One can see this as an independent
increment property, where the ``increment'' is the conformal equivalence $(g_s-W_s)$.

For any $\kappa>0$, there exists a continuous process $\gamma$ taking values in
$\overline\H$ that generates the hulls $(K_t)$ of $\SLE(\kappa)$ in
the following sense: a.s., $\H\setminus K_t$ is the unbounded
connected component of $\H\setminus\gamma_{[0,t]}$ for any $t\geq 0$;
this process $\gamma$ is the trace of the $\SLE$ (see \cite{RS01},
\cite{LSW2} for the case $\kappa=8$).

If $\kappa\leq 4$, the trace is a.s. simple; this is no longer the
case if $\kappa>4$ (see \cite{RS01}). Hence, if $\kappa>4$, the outer
boundary of an $\SLE(\kappa)$ hull $K_t$ is strictly included in the
image of the trace $\gamma_{[0,t]}$. Furthermore, if $4<\kappa<8$, an
$\SLE(\kappa)$ hull $K_t$ has cutpoints with positive probability (\cite{Bef}).

Other probability measures on processes can be translated in laws on
increasing families of hulls by means of Loewner equations. For
$n\in\N$ and a family of parameters $\underline\rho=(\rho_1,\dots, \rho_n)$, consider the SDEs
$$\left\{\begin{array}{ccl}
dW_t&=&\sqrt\kappa dB_t+\sum_{i=1}\frac{\rho_idt}{W_t-Z^{(i)}_t}\\
dZ^{(i)}_t&=&\frac {2dt}{Z^{(i)}_t-W_t}
\end{array}\right.
$$
for $i=1\dots n$, where $B$ is a standard Brownian motion, with initial conditions
$W_0=0, Z^{(i)}_0=x_i\in\R$. Then the image of the process $W$ under
the Loewner equations define $\SLE(\kappa,\underline\rho)$. The
$\SLE(\kappa,\rho)$ processes were introduced in \cite{LSW3}; see also
\cite{W2,D4}. In fact, we will use only the cases $n=1,2$; questions
of definiteness will be discussed when needed. 

\subsection{Properties of $\SLE(\kappa,\kappa-4)$ and
  $\SLE(\kappa,\kappa-4,\kappa-4)$ processes}

In this section, we collect some of the properties of  $\SLE(\kappa,\kappa-4)$ and
  $\SLE(\kappa,\kappa-4,\kappa-4)$ that we will use later. Let us
  stress that almost all of these properties seem to have no direct equivalent
  for other choices of the $\rho$ parameter.
% For general background on
%  $\SLE$ and $\SLE(\kappa,\rho)$, see Sections 5.2, 5.5, 5.6, and
%  references therein. 
Let $\kappa>4$ be fixed, and $d=1-\frac
  4\kappa$. Recall that a $\Beta(a,b)$ law is a probability law
  supported on $[0,1]$ with density:
$$B(a,b)^{-1}\ind_{x\in (0,1)}x^{a-1}(1-x)^{b-1}$$
where $B(a,b)=\Gamma(a)\Gamma(b)/\Gamma(a+b)$, and $a$, $b$ are positive.

\begin{Prop}\label{BP1}
Let $(W_t,O_t)_t$ be the driving process of an $\SLE(\kappa,\kappa-4)$
starting from $(0,y)$, $y>0$ and $({\mc F}_t)$ the associated natural
filtration. For $x\in (0,y)$, let $\tau_x\in
(0,\infty]$ be the swallowing time of $x$. Then:\\ 
(i) $\P(\tau_x=\infty)=(x/y)^d$.\\
(ii) The law of the $\SLE$ conditionally on $\{\tau_x=\infty\}$ is that
of an $\SLE(\kappa,\kappa-4)$ starting from $(0,x)$.\\
(iii) Let $Z$ be the rightmost swallowed point of $(0,y)$. Then $Z/y$ is
an ${\mc F}_\infty$-measurable sample from the $\Beta(d,1)$ law.\\
(iv) Let $\tau_Z$ be the a.s. finite ${\mc F}_\infty$-measurable
random time at which $Z$ is swallowed. For $z\in (0,1)$, let $Q(z,.)$
denote the probability measure on bivariate $(W_t,O_t)$ processes
obtained as the concatenation of an $\SLE(\kappa,-4)$ starting from
$(0,z)$, which is defined up to time $\tau_z$, and an independent
$\SLE(\kappa,\kappa-4)$ starting from
$(W_{\tau_z},W_{\tau_z}^+)$. Then $\omega\mapsto Q(Z(\omega),.)$
defines a regular conditional probability of the original
$\SLE_{(0,y)}(\kappa,\kappa-4)$ process w.r.t. $\sigma(Z)$.
\end{Prop}  
\begin{proof}
(i) It is immediate to check that the semimartingale:
$$t\mapsto \left(\frac{g_t(x)-W_t}{O_t-W_t}\right)^d,$$
which is defined up to $\tau_x$, is a local martingale taking values
in $[0,1]$, hence is a martingale. If $\tau_x<\infty$, as $y$ is never
swallowed, it appears that as $t\nearrow\tau_x$, this martingale goes
to zero (as the extremal distance between $g_t^{-1}((W_t,g_t(x)))$
and $(y,\infty)$ in $\H\backslash K_t$ goes to
infinity). Conversely, if $\tau_x=\infty$, this extremal distance
goes to 0 (by transience of the $\SLE$), so the limit of the
martingale is 1. One concludes with the optional stopping theorem (as
e.g. in \cite{W1}, 3.2).

(ii) See \cite{D4}, Section 5.

(iii) The density follows from (i), since, from the definition of $Z$, one has
: $\{Z<x\}=\{\tau_x=\infty\}$.

(iv)  We are interested in the law of the $\SLE$ conditioned on
$Z$. Suppose $Z\in(z-\eps,z)$ for some $z\in (0,1]$; then $z$ is not
swallowed, so from (ii), an $\SLE_{(0,y)}(\kappa,\kappa-4)$ conditionally on
$Z\in (z-\eps,z)$ has the law of an $\SLE_{(0,z)}(\kappa,\kappa-4)$ conditionally on $Z>z-\eps$.
We argue as in \cite{D4}, Section 5; if ${\bf P}$ denotes the
original probability measure (corresponding to a
$\SLE(\kappa,\kappa-4)$ starting from $(0,y)$\ ) and ${\bf Q}={\bf
  P}(.|Z\in(z-\eps,z))$, then there exists a process $\tilde B$, which
is a standard Brownian motion under ${\bf Q}$, such that:
%\begin{align}\label{sdeps}
%dW_t=&\sqrt\kappa d\tilde B_t+\frac{\kappa-4}{g_t(z)-W_t}dt\nonumber\\
%&-(\kappa-4)\frac{g_t(z)-g_t(z-\eps)}{(g_t(z)-W_t)^2}U_t^{-4/\kappa}\left(1-U_t^{1-4/\kappa}\right)^{-1}dt
%\end{align}
\begin{equation}\label{sdeps}
dW_t=\sqrt\kappa d\tilde B_t+\frac{\kappa-4}{g_t(z)-W_t}dt
-(\kappa-4)\frac{g_t(z)-g_t(z-\eps)}{(g_t(z)-W_t)^2}U_t^{-4/\kappa}\left(1-U_t^{1-4/\kappa}\right)^{-1}dt
\end{equation}
where $U_t=(g_t(z-\eps)-W_t)/(g_t(z)-W_t)$. As $\eps\searrow 0$, the
drift term in this equation tends to $-4/(g_t(z)-W_t)dt$. This limit
SDE defines an $\SLE_{(0,z)}(\kappa,-4)$ process; this indicates that
an  $\SLE_{(0,y)}(\kappa,\kappa-4)$ conditioned by $Z\in dz$, and
stopped on swallowing $Z$ at time $\tau_Z$, is an
$\SLE_{(0,z)}(\kappa,-4)$ stopped at time $\tau_z$. Note that this
last process is well defined up to time $\tau_z$, and that in this
case $(g_t(z)-W_t)/\sqrt\kappa$ is a standard $\BES(d)$ process stopped on hitting 0.
%Let $\eta$
%be an arbitrary small number, and $\tau_{\eta,\eps}=\min(\inf(t\geq
%0,g_t(z-\eps)-W_t\leq\eta),1/\eta)$. If $(W_t)$ is the driving
%process of a $\SLE_{(0,y)}(\kappa,\kappa-4)$ conditioned on $Z\in
%(z-\eps,z)$, then, as $\eps\searrow 0$, $(W_t)_{t\leq\tau_{\eta,\eps}}$ converges in law to
%the driving process of a $\SLE_{(0,z)}(\kappa,-4)$ stopped at time
%$\tau_{\eta,0}$; in this case, $(g_t(z)-W_t)/\sqrt\kappa$ is a
%standard $\BES(d)$ process stopped on hitting 0. 

We make the previous limiting argument a bit more precise. Let $\eps_0>0$ be
fixed. Consider a process $R=(R_1,R_2,R_3)$ taking values in the Banach space
$C_0([0,\eps_0],\R^3)$ (for uniform convergence w.r.t. some norm of $\R^3$), that satisfies the following SDE:
$$dR_t=(\sqrt\kappa,0,0)d\tilde B_t+\left(b_1(R_t),\frac
2{R_{2,t}-R_{1,t}},\frac 2{R_{3,t}-R_{1,t}}\right)dt$$
with initial conditions $R_0(\eps)=(0,z-\eps,z)$. The drift term
$b_1(r)$ is defined by:
$$b_1(r)=\frac{\kappa-4}{r_3-r_1}-(\kappa-4)\frac{r_3-r_2}{(r_3-r_1)^2}\left(\frac{r_2-r_1}{r_3-r_1}\right)^{-4/\kappa}\left(1-\left(\frac{r_2-r_1}{r_3-r_1}\right)^{1-4/\kappa}\right)^{-1}$$
with $b_1(r)(0)=\lim_{\eps\searrow 0} b_1(r)(\eps)=-4/(r_3-r_1)$. Now
let $\eta>0$ be an arbitrary small number, and let $\Delta_\eta$ be
the following closed subset of $C_0([0,\eps_0],\R^3)$:
$$\Delta_\eta=\{(r_1,r_2,r_3): r_1\leq r_2\leq r_3,
r_3(\eps)-r_2(\eps)\leq \eps, r_2-r_1\geq \eta\}.$$
Let $b_2(r)=2/(r_2-r_1)$, $b_3(r)=2/(r_3-r_1)$, $b=(b_1,b_2,b_3)$. It
is obvious that $b_2$ and $b_3$ define uniformly Lipschitz continuous
functions from $\Delta_\eta$ to $C_0([0,\eps_0],\R)$. As for $b_1$, note
that $r\mapsto 1/(r_3-r_1)$ and $r\mapsto (r_3-r_2)/(r_3-r_1)$ are
uniformly Lipschitz continuous on $\Delta_\eta$, and that for all
$\eps\in [0,\eps_0]$, $r\in\Delta_\eta$\ :
$$0\leq\frac {r_3-r_2}{r_3-r_1}=1-\frac {r_2-r_1}{r_3-r_1}\leq
1-\frac{r_2-r_1}{r_2-r_1+\eps}\leq 1-\frac 1{1+\eps/\eta}\leq \frac{\eps_0}{\eps_0+\eta}$$
Since the function
$$u\mapsto u(1-u)^{-4/\kappa}\left(1-(1-u)^{1-4/\kappa}\right)^{-1}$$
extends to a $C^1$ function on $[0,1[$, it defines a uniformly
Lipschitz continuous function on $[0,\eps_0/(\eps_0+\eta)]$. It
follows that $b: \Delta_\eta\rightarrow C_0([0,\eps_0],\R^3)$ is
uniformly Lipschitz continuous. So existence and strong unicity hold
for the SDE
$$dR_t=(\sqrt\kappa,0,0)d\tilde B_t+b(R_t)dt,$$
 up to the first exit of $\Delta_\eta$. Since 
$$d(R_{3,t}-R_{2,t})=-\frac {2(R_{3,t}-R_{2,t})dt}{(R_{3,t}-R_{1,t})(R_{3,t}-R_{1,t})}$$
the difference $(R_3-R_2)$ is decreasing, and a solution of the SDE is
well defined until $(R_2-R_1)$ reaches $\eta$. Letting $\eta$ go to
$0$, one gets existence and strong unicity up to explosion of $(R_2-R_1)^{-1}$.

Now, for any $\eps\in (0,\eps_0]$, $(R_t(\eps))_t$ has the same law as
$(W_t,g_t(z-\eps),g_t(z))_t$, where $(W_t)$ is the driving process of an
$\SLE_{(0,y)}(\kappa,\kappa-4)$ conditionally on $Z\in (z-\eps,z)$,
and $(g_t)$ are the corresponding conformal equivalences; for
$\eps=0$, one gets an $\SLE_{(0,z)}(\kappa,-4)$. So the driving
process of an $\SLE_{(0,y)}(\kappa,\kappa-4)$ conditionally on $Z\in
(z-\eps,z)$ and stopped at time $\tau_{z-\eps}$  converges in law as
$\eps\searrow 0$ to that of an $\SLE_{(0,z)}(\kappa,-4)$ stopped as
time $\tau_z$ (after which it is no longer defined).
 
The future after $\tau_z$ is easier to handle. Indeed, from the Strong
Markov property, an $\SLE_{(0,y)}(\kappa,\kappa-4)$ can be described
as the concatenation of this process stopped at time $\tau_{z-\eps}$
and an
$\SLE_{(W_{\tau_{z-\eps}},g_{\tau_{z-\eps}}(y))}(\kappa,\kappa-4)$,
  independent of the former conditionally on its starting
  state. Conditioning by $Z\in (z-\eps,z)$, one gets the
  concatenation of the original conditional process stopped at time
  $\tau_{z-\eps}$ and an $\SLE(\kappa,\kappa-4)$ starting from
  $(W_{\tau_{z-\eps}},g_{\tau_{z-\eps}}(z))$ (as follows from (ii)),
    independent of the former conditionally on this starting state;
    note that $W_{\tau_{z-\eps}}\leq g_{\tau_{z-\eps}}(z)\leq
      W_{\tau_{z-\eps}}+\eps$.

As $\eps\searrow 0$, the conditional process up to time
$\tau_{z-\eps}$ converges to an $\SLE_{(0,z)}(\kappa,-4)$ stopped at time
$\tau_z$, while the process after $\tau_{z-\eps}$ converges in law to
an $\SLE_{(W_{\tau_z},W_{\tau_z}^+)}(\kappa,\kappa-4)$, independent of
the former conditionally on $W_{\tau_z}$. This concludes the proof.
\end{proof}

We have mentioned that the $\SLE_{(0,y)}(\kappa,\kappa-4)$ could be
seen as a standard $\SLE(\kappa)$ conditioned ``never to swallow $y$''
(see \cite{D4}, Section 5). Likewise, a
$\SLE_{(0,y_1,y_2)}(\kappa,\kappa-4,\kappa-4)$, $y_1<0<y_2$, can be
interpreted as an $\SLE(\kappa)$ conditioned ``never to swallow $y_1$
or $y_2$''. On a heuristic level, this explains why
$\SLE(\kappa,\kappa-4,\kappa-4)$ should have properties similar to
those previously described for $\SLE(\kappa,\kappa-4)$.

\begin{Prop}\label{BP2}
Let $(W_t,O^{1}_t,O^{2}_t)_t$ be the driving process of an $\SLE(\kappa,\kappa-4,\kappa-4)$
starting from $(0,y_1,y_2)$, $y_1<0<y_2$; let $({\mc F}_t)$ the associated natural
filtration. For some $x_1\in [y_1,0)$, $x_2\in (0,y_2]$, let $\tau_i\in
(0,\infty]$ be the swallowing time of $x_i$, $i=1,2$.\\ 
(i) The following formula holds:
$$\P(\tau_1=\infty,\tau_2=\infty)=\left(\frac{x_1x_2}{y_1y_2}\right)^d\left(\frac{x_2-x_1}{y_2-y_1}\right)^{\frac
  \kappa 2
  d^2}.$$
(ii) The law of the $\SLE$ conditionally on $\{\tau_1=\infty,\tau_2=\infty\}$ is that
of an $\SLE(\kappa,\kappa-4,\kappa-4)$ starting from $(0,x_1,x_2)$.\\
(iii) Let $Z_1$ (resp. $Z_2$) be the leftmost (resp. rightmost)
swallowed point of $(y_1,0)$ (resp. $(0,y_2)$). Then $(Z_1,Z_2)$ is
an ${\mc F}_\infty$-measurable r.v. with density:
$$\ind_{y_1\leq x_1<0<x_2\leq y_2} d^2(\frac\kappa 2-1)\left(-\frac 1{x_1x_2}+\frac{\kappa/2-4}{(x_2-x_1)^2}\right)\left(\frac{x_1x_2}{y_1y_2}\right)^d\left(\frac{x_2-x_1}{y_2-y_1}\right)^{\frac
  \kappa 2
  d^2}dx_1dx_2.$$
\end{Prop}
\begin{proof}
(i) As in the previous proposition, there is only to check, using
It\^o's formula, that the following semimartingale
$$\left(\frac{(g_t(x_1)-W_t)(g_t(x_2)-W_t)}{(g_t(y_1)-W_t)(g_t(y_2)-W_t)}\right)^d\left(\frac{g_t(x_2)-g_t(x_1)}{g_t(y_2)-g_t(y_1)}\right)^{\frac\kappa
  2d^2}$$
is a local martingale.

(ii) Using (i) and Girsanov's theorem, one can give a proof that
follows exactly that of \ref{BP1} (ii). Note that if $(M_t)$ is
the martingale considered in (i), and $(B_t)$ is the driving Brownian
motion of the $\SLE$, then:
$$\frac {d\langle M_t,\sqrt\kappa
  B_t\rangle}{M_t}=(\kappa-4)\left(\left(\frac 1{g_t(x_1)-W_t}+\frac 1{g_t(x_2)-W_t}\right)-\left(\frac 1{g_t(y_1)-W_t}+\frac 1{g_t(y_2)-W_t}\right)\right)dt
$$
which is the difference between the drift term of an
$\SLE_{(0,x_1,x_2)}(\kappa,\kappa-4,\kappa-4)$ and that of an
$\SLE_{(0,y_1,y_2)}(\kappa,\kappa-4,\kappa-4)$. We shall give a
clean interpretation of this identity later.

(iii) This follows directly from (i), since
$\{x_1<Z_1<0<Z_2<x_2\}=\{\tau_1=\tau_2=\infty\}$, from the definition
of $Z_1,Z_2$.

\end{proof}

\section{Frontier points of $\SLE(\kappa,\kappa-4)$ processes}

As before, $\kappa>4$ is fixed, and $d=1-\frac 4\kappa$. 
Let $(W_t,O_t)_{t\geq 0}$ be the driving process of an
$\SLE(\kappa,\kappa-4)$ starting from $(0,0^+)$; $(K_t)$, $(g_t)$, $(\gamma_t)$ are
as usual; $({\mc F}_t)$ is the natural filtration of the Brownian
motion driving the $\SLE$. Then, if
$Y_t=O_t-W_t$, $Y/\sqrt\kappa$ is a standard Bessel process of
dimension $(2+d)$ starting from $0$. Indeed, it satisfies the SDE:
$$dY_t=dB_t+\frac{d+1}2dt$$
where $B$ is the (standard) Brownian motion driving the $\SLE$. For
general background on Bessel processes, see \cite{RY}, XI.1.

Let $t$ be a positive time, and $\sigma_t$ be the last time
before $t$ when the tip of hull $\gamma_t$ was on the right-boundary of the whole
hull $K_\infty$; this is obviously not a stopping time, and is
analogous to the future infimum of a transient Bessel process for instance. Now define:
$$X_t=g_t(\gamma_{\sigma_t})-W_t.$$
From the Markov property of $\SLE$, $(g_t(K_{t+s}\backslash
K_t)-W_t)_s$ defines an $\SLE_{(0,Y_t)}(\kappa,\kappa-4)$,
independent from ${\mc F}_t$ conditionally on $Y_t$. Then $X_t$
is the rightmost point on $(0,Y_t)$ swallowed by this $\SLE$. It
follows from scale invariance that $X_t/Y_t$ is independent from
$\mc{F}_t$, and we have seen in Proposition \ref{BP1} that
$X_t/Y_t\stackrel{\mc L}{=} \Beta(d,1)$.

As $X_t=Y_t(X_t/Y_t)$, these
two variables being independent, one may work out the law of
$X_t$. Indeed, for a standard Bessel process of dimension $\delta$
starting from 0,
$$\P(\BES(\delta)_t\in dr)=\frac{2r^{\delta-1}e^{-r^2/2t}dr}{\Gamma(\delta/2)(2t)^{\delta/2}}.$$
It follows that:
\begin{align*}
\P(X_t/\sqrt\kappa\in dx)&=\int_x^\infty \P(X_t/\sqrt\kappa\in dx|Y_t/\sqrt\kappa\in
dy)\P(Y_t/\sqrt\kappa\in dy)\\
&=\int_x^\infty
d\frac{x^{d-1}dx}{y^d}.\frac{2y^{d+1}e^{-y^2/2t}}{\Gamma(1+d/2)(2t)^{1+d/2}}dy\\
&=\frac{dx^{d-1}e^{-x^2/2t}}{\Gamma(d/2+1)(2t)^{d/2}}dx=\frac{2x^{d-1}e^{-x^2/2t}}{\Gamma(d/2)(2t)^{d/2}}dx.
\end{align*}

So it turns out that $X_t/\sqrt\kappa$ is distributed like a standard
$\BES(d)$ starting from $0$ and taken at time $t$; taking squares, one
may see this as a special case of beta-gamma algebra (see
\cite{CPY}). Let $({\mc X}_t)$ be the filtration generated by this
process. It appears readily that ${\mc X}_t\subset {\mc F}_\infty$,
and ${\mc X}_t\nsubseteq {\mc F}_s$ for any $s\geq 0$.

The computation of the distribution of $X_t$ at a fixed time $t$ suggests the following generalization:

\begin{Prop}\label{BP3} The process $(X_t/\sqrt\kappa)_{t\geq 0}$ is a standard $\BES(d)$ process.
\end{Prop}
\begin{proof}
%The proof relies on regeneration for $\SLE(\kappa,\kappa-4)$, as well
%as the law of an $\SLE(\kappa,\kappa-4)$ starting from $(0,1)$
%conditionally on $Z\in(0,1)$, the rightmost swallowed point (which has
%distribution $\Beta(d,1)$). 
Let $(Q_t)$ designate the semigroup of a standard $d$-dimensional Bessel process. We
will prove that for all $0\leq s\leq t$, the law of $X_t/\sqrt\kappa$
conditionally on ${\mc X}_s$ is $Q_{t-s}(X_s/\sqrt\kappa,.)$. The
computation above shows that this holds for $0=s\leq t$. 

Suppose now that $0<s<t$ are fixed. Almost surely, $X_s$ is
positive. Conditionally on ${\mc X}_s\vee{\mc Y}_s$, the translated
process $(Y_{s+u})_u$ corresponds to an
$\SLE_{(0,Y_s)}(\kappa,\kappa-4)$ conditioned on its rightmost
swallowed point being $X_s$. Indeed, ${\mc X}_s\vee{\mc
  Y}_s=\sigma({\mc Y}_s,X_s)$, since for $r\leq s$,
$X_r=g_r(\gamma_{\sigma_r})-W_r$, and: 
$$\sigma_r=\sup\left(u:\ u\leq r,\ \gamma_u\in
g_s^{-1}((\gamma_{\sigma_s},\infty))\right).$$ 
From Proposition \ref{BP1} (iv), we know that $(X_{s+u})_u$
(stopped on hitting 0) corresponds to an $\SLE_{(0,X_s)}(\kappa,-4)$
independent of ${\mc X}_s\vee{\mc Y}_s$ conditionally on
$(Y_s,X_s)$. But $Y_s$ does not appear any longer in the conditional
law (which essentially follows from Proposition \ref{BP1} (ii)\,), so
we have proved that the law $(X_{s+u}/\sqrt\kappa)_u$ stopped on hitting 0
conditionally on ${\mc X}_s\vee{\mc Y}_s$ is that of a standard
$\BES(d)$ process starting from $X_s/\sqrt\kappa$, and stopped on hitting $0$.

Consider now $T$, an ${\mc X}\vee {\mc Y}$-stopping time a.s. supported on
zeroes of $X$. We shall prove that $(g_{T+u}\circ
g_T^{-1}(W_T^+)-W_{T+u},X_{T+u})_u$ is a copy of $(Y_u,X_u)_u$ independent
from $({\mc X}\vee {\mc Y})_T$. Indeed, let $0<t_1<t_2$; on the event
$T\in (t_1,t_2)$, $T$ is $({\mc X}\vee {\mc
  Y})_{t_2}=\sigma(X_{t_2},{\mc Y}_{t_2})$-measurable. The process
$(g_{t_2+u}\circ g_{t_2}^{-1}(X_{t_2}^+ +W_{t_2})-W_{t_2+u},X_{t_2+u})_u$ is
distributed as the process $(Y_u,X_u)$ for an
$\SLE_{(0,X_{t_2})}(\kappa,\kappa-4)$ conditioned on its right-most
swallowed point being $X_{t_2}$, as described in Proposition \ref{BP1}
(iv); this translated process is independent from $({\mc X}\vee {\mc
  Y})_{t_2}$ conditionally on $X_{t_2}$. As $t_1\nearrow t_2$, since
$X_T=0$ and $X$ is continuous, the law of the translated process
converges to that of a copy of $(Y_u,X_u)_u$ independent from $({\mc
  X}\vee {\mc Y})_{t_2}$. As this holds for any positive $t_2$, the
claim follows.

There only remains to put things together. Let $0\leq s\leq t$, and
consider the ${\mc X}\vee {\mc Y}$-stopping time $T=\inf(u\in [s,t],
X_u=0)$, with the convention $\inf\varnothing=\infty$. Let $(\tilde
X_u/\sqrt\kappa)_{s\leq u\leq t}$ be a standard $\BES(d)$ process
starting from $\tilde X_s=X_s$, and $\tilde T=\inf(u\in [s,t],
\tilde X_u=0)$. Then the $[s,t]\cup\{\infty\}$-valued variables $T$ and
$\tilde T$ have the same distribution. Either $T=\infty$, and
$X_t$ is distributed as $\tilde X_t$ conditionally on
$\tilde T=\infty$, or $T\leq t$; in the latter case, we have seen
that $X_t/\sqrt\kappa$ is independent from $\sigma((X_u)_{s\leq
u\leq T})$ conditionally on $T$, and has conditional distribution
$Q_{t-T}(0,.)$. From the Strong Markov property for Bessel processes, 
the following disintegration holds:
$$Q_{t-s}(x,.)=\int_s^t Q_{t-u}(0,.)\P(\tilde T\in
du)+Q_{t-s}(x,.|\tilde T=\infty)\P(\tilde T=\infty)$$
where $Q_{t-s}(x,.|\tilde T=\infty)$ designates the law of $\tilde
X_t$ conditionally on $\inf_{s\leq u\leq t}(\tilde X_u)>0$. Hence we have proved that (a regular version of) the
distribution of $X_t/\sqrt\kappa$ conditionally on $({\mc X}\vee {\mc
  Y})_s$ is $Q_{t-s}(X_s/\sqrt\kappa,.)$, so that $(X_t/\sqrt\kappa)_t$ is
a standard $d$-dimensional Bessel process.
\end{proof}

We further discuss the relationship between the processes $X$ and
$Y$. These two processes are Bessel processes of dimension $d$ and
$(d+2)$ in their natural filtration; but these two filtrations differ,
since ${\mc X_t}\nsubseteq {\mc Y}_t$. We will see later that in fact
${\mc Y}_t\subset{\mc X}_t$. Another important feature is the
intertwining of these two Markov processes. Indeed, we have mentioned that
$(X_t/Y_t)$ is distributed as a $\Beta(d,1)$ variable, independently of  
${\mc Y}_t$. So if $(P^X_t)$ and $(P^Y_t)$ denote the Markov
semigroups of $X$ and $Y$ respectively (starting from 0), and
$\Lambda$ is the Markov transition kernel specified by its action on
bounded Borel functions on $\R^+$:
$$(\Lambda f)(y)=d\int_0^1 u^{d-1}f(uy)\,du$$ 
then the following intertwining relation holds: $P^Y\Lambda=\Lambda
P^X$. Indeed, we have:
\begin{align*}
\E(f(X_t)|{\mc Y}_s)&=\E(f(X_t)|{\mc X}_s\vee{\mc Y}_s|{\mc
  Y}_s)=\E(P^X_{t-s}f(X_s)|{\mc Y}_s)=(\Lambda P^Xf)(Y_s)\\
&=\E(f(X_t)|{\mc Y}_t|{\mc
  Y}_s)=\E((\Lambda f)(Y_t)|{\mc
  Y}_s)=(P^Y\Lambda f)(Y_s)
\end{align*}
as in \cite{CPY}, Proposition 2.1.

The definition and study of the process $X$ have been motivated by
questions related to the final right boundary of the $\SLE(\kappa,\kappa-4)$
process. Notice that this boundary grows only when $X$ vanishes, so it is now
natural to turn to the excursion theory of $X$ (away from 0). In particular the
local time at 0 of $(X_t)$ will provide an adequate measure of the
``size'' of the right boundary. 

Recall that for $(\rho_t)$ a Bessel process of
dimension $\delta$, there exists a bicontinuous family of local times
$(L^a_t)$ such that for any bounded Borel function $f$,
$$\int_0^t f(\rho_u)du=\int_0^\infty f(a)L^a_ta^{\delta-1}da.$$
(For Bessel local times, see e.g. \cite{RY}, XI.1.25). Since $\rho$ has Brownian scaling, one gets: 
$$(L^{\lambda
  a}_{\lambda^2 t})\stackrel{{\mc L}}=\lambda^{2-\delta}(L^a_t)$$
as a bicontinuous process of the two variables $t$ and $a$. As a consequence, if $(L^0_t)_{t\geq 0}$ is the local time at $0$ of
$(X_t)$, and if $\tau_.$ is its right-continuous inverse, then $\tau.$
is a stable subordinator with index $\nu=1-\frac d2=\frac 12+\frac
2\kappa$. 

From the definition of $X$, it appears that the set of (global) frontier times
of $\SLE(\kappa,\kappa-4)$ (i.e. times at which the trace lies on the
right-boundary of the total hull $K_\infty$) is the zero set of
$X$. Note that the a.s. Hausdorff dimension of frontier times has
been derived in \cite{Bef}. Recall that the $\varphi$-measure of a
set $Z$ is defined as:
$$\varphi-m(Z)=\liminf_{\eps>0}\left(\sum_1^\infty\varphi(\diam
(J_i)), Z\subset\bigcup_1^\infty J_i, \diam (J_i)<\eps\right).$$

\begin{Cor}
The set of frontier times of $\SLE(\kappa,\kappa-4)$ has
a.s. Hausdorff dimension $(1/2+2/\kappa)$. More precisely, the
$\varphi_\kappa$-measure of frontier times is a.s. positive and
locally finite, where:
$$\varphi_\kappa(\eps)=\eps^{1/2+2/\kappa}(\log|\log\eps|)^{1/2-2/\kappa}.$$
Moreover, the $\varphi_\kappa$-measure of frontier times up to time
$t$ is a multiple of the local time of $X$ at $0$.
\end{Cor}
\begin{proof}
The inverse of the local time of $X$ at $0$ is a stable subordinator
with index $(1/2+2/\kappa)$. So the zero set of $X$ has the same distribution
as the zero set of a stable L\'evy process with index
$(1/2-2/\kappa)^{-1}$. Hence one can apply the results of \cite{TW}.
\end{proof}

We will also need some facts about principal values for $\delta\in
(0,1)$, that we presently recall (see \cite{BerComp}). It is classical that one can consider a version of $(L^a_t)$ that is bicontinuous and
H\"older in the space variable. More precisely, one can use the
Biane-Yor identity (see \cite{RY}) to translate properties of Brownian
local times (which are $(1/2-\eps)$-H\"older in the space variable for
any $\eps>0$) into results on Bessel local times. Hence one can define:
$$p.v.\int_0^t\frac{ds}{\rho_s}\stackrel{\text{def}}{=}\int_0^\infty (L^a_t-L^0_t)a^{\delta-2}da.$$
As a function of $t$, this is continuous, and $C^1$ on $\{s:\rho_s\neq
0\}$, with derivative $1/\rho_s$; though, it is not monotonic.

For $\delta>1$, $\rho$ is a semimartingale with decomposition:
$$\rho_t=B_t+\frac{\delta-1}2\int_0^t\frac{ds}{\rho_s}$$
where $B$ is a standard Brownian motion. If $0<\delta<1$, $\rho$ is no
longer a semimartingale, but it still is a
Dirichlet process (i.e. the sum of a square-integrable local martingale and a process
with zero quadratic variation):
$$\rho_t=B_t+\frac{\delta-1}2p.v.\int_0^t\frac{ds}{\rho_s}$$
where $B$ is a standard Brownian motion and the principal value has
zero quadratic variation.

We turn back to the previous situation: $(W_t,O_t)$ is the driving
mechanism of an $\SLE_{(0,0^+)}(\kappa,\kappa-4)$, $d=1-4/\kappa$,
$Y_t=O_t-W_t$, $X_t=g_t(\gamma_{\sigma_t})-W_t$.

\begin{Prop} The following identity holds a.s.:
$$W_t=-X_t+p.v.\int_0^t \frac{2ds}{X_s}=-Y_t+\int_0^t\frac{2ds}{Y_s}.$$
Consequently, ${\mc Y}_t={\mc F}_t\subset {\mc X}_t$ for all $t\geq0$.
\end{Prop}
\begin{proof}
The second half of the identity is just the definition of
$\SLE(\kappa,\kappa-4)$.
A consequence of Loewner's equation and of the definition of $(X_t)$ is that:
$$S_t=X_t+W_t-p.v.\int_0^t \frac{2ds}{X_s}$$
is a continuous process that is constant on the excursions of $X$
(recall that $X/\sqrt\kappa$ is a standard $d$-dimensional Bessel
process). Indeed, $X_t=g_t(\gamma_{\sigma_t})-W_t$, and $\sigma_.$ is
constant on the excursions of $X$. Define $\tilde S_u=S_{\tau_u}$ where $\tau_u$ is the right-continuous
inverse of $(L^0_t)$, the local time at 0 of $X$.
Since $(S_t)$  is continuous and constant on the excursions of $X$,
the process $\tilde S$ has a continuous version. Indeed, $\tilde S$ is
c\`adl\`ag; suppose that is has a discontinuity at $u$. Necessarily,
$\tau_u>\tau_u^-$, and $(\tau_u^-,\tau_u)$ is an excursion interval for
$X$, so $S_{\tau_u^-}=S_{\tau_u}$. Since the local time at 0 of $X$
is instantaneously increasing at the right of $\tau_u$ and at the left of
$\tau_u^-$, and $S$ is continuous at $\tau_u$ and $\tau_u^-$, so is
$\tilde S$.

Moreover, for $u>0$, $\tau_u$ is a stopping time for $X$ such that
$X_{\tau_u}=0$ a.s.\ . As we have seen in the proof of Proposition \ref{BP3}, this implies that
$(X_{\tau_u+t},W_{\tau_u+t}-W_{\tau_u})_t$ has the same law as
$(X_t,W_t)_t$ and is independent from ${\mc X}_{\tau_u}\vee {\mc F}_{\tau_u}$. From the definition of $S$ and $\tilde S$, it appears that $\tilde S_{u+v}=\tilde S_u+\tilde S'_v$ where
$\tilde S'$ is a copy of $\tilde S$ independent from $\sigma((X_s,W_s)_{s\leq
  \tau_u})$. So $\tilde S$ is a continuous L\'evy process, i.e. a
Brownian motion with drift. But $S$ has Brownian scaling, and $\tau_.$
has index $\nu=1-d/2$, so $\tilde S$ has index $2\nu\in (1,2)$. Hence
$\tilde S$ and $S$ vanish identically, which proves the first half
of the identity.

Furthermore, $W_t$ is ${\mc X}_t$ measurable for all $t$, so  ${\mc
  Y}_t={\mc F}_t\subset {\mc X}_t$ for all $t\geq0$. 
\end{proof}

Before proceeding with the study of $\SLE(\kappa,\kappa-4)$, we sum up some of these results in an $\SLE$-free formulation.

\begin{Cor}
For $d\in (0,1)$, let $\rho_1$ (resp. $\rho_2$) be a standard Bessel
process of dimension $d$ (resp. $(d+2)$), and ${\mc F}^1$ (resp. ${\mc
F}^2$) be its natural filtration. Then there exists a coupling
of $\rho_1$ and $\rho_2$ such that:\\
(i) ${\mc F}^2_t\subset {\mc F}^1_t$ for all $t\geq 0$.\\
(ii) The process $(\rho_2-\rho_1)$ has zero quadratic variation.\\
(iii) The semigroups of $\rho_1,\rho_2$ are intertwined by
multiplication by a $\Beta(d,1)$ law. More precisely, if $f$ is some
bounded Borel function on $\R^+$, then:
$$\E(f(\rho_{1,t})|{\mc F}^2_t)=\int_0^1f(u\rho_{2,t})d\,u^{d-1}\,du.$$
(iv) The following identity holds a.s. for all $t\geq 0$:
$$\rho_{1,t}+\frac{d-1}2 p.v.\int_0^t\frac {ds}{\rho_{1,s}}=\rho_{2,t}+\frac{d-1}2 \int_0^t\frac {ds}{\rho_{2,s}}.$$
(v) Let $\rho_{1,t}=B^1_t+\frac{d-1}2p.v.\int_0^t\frac
{ds}{\rho_{1,s}}$ and $\rho_{2,t}=B^2_t+\frac{d+1}2\int_0^t\frac
{ds}{\rho_{2,s}}$ be the Dirichlet process decompositions of $\rho_1$
and $\rho_2$. Then:
\begin{align*}
B^2_t&=\E\left(B^1_t|{\mc F}^2_t\right)\\
d\int_0^t\frac{ds}{\rho_{2,s}}&=(d-1)\E\left(\left.p.v.\int_0^t\frac{ds}{\rho_{1,s}}\right|{\mc F}^2_t\right).
\end{align*}
\end{Cor}
\begin{proof}
With $d=1-4/\kappa$, define $\rho_1=X/\sqrt\kappa$ and
$\rho_2=Y/\sqrt\kappa$, where $X$ and $Y$ are as above. We have
already proved (i), (iii), (iv). For (ii), note that:
$$\sqrt\kappa(\rho_{2,t}-\rho_{1,t})=Y_t-X_t=p.v.\int_0^t\frac {2ds}{X_s}-\int_0^t\frac {2ds}{Y_s}.$$
(v) Taking conditional expectations in (iv), one gets:
$$\E\left(\left.\rho_{1,t}+\frac{d-1}2 p.v.\int_0^t\frac
  {ds}{\rho_{1,s}}\right|{\mc F}^2_t\right)=\rho_{2,t}+\frac{d-1}2 \int_0^t\frac {ds}{\rho_{2,s}}$$
since the right-hand side is ${\mc F}^2_t$-measurable. It follows
that:
$$\E(B^1_t|{\mc
  F}^2_t)-B^2_t=d\int_0^t\frac{ds}{\rho_{2,s}}-(d-1)\E\left(\left.p.v.\int_0^t\frac{ds}{\rho_{1,s}}\right|{\mc F}^2_t\right).$$
The left-hand side is a continuous square-integrable ${\mc F}^2$-
martingale, while the right-hand side has zero quadratic variation; so
both sides vanish identically.
\end{proof}

Taking expectations in (v), one gets:
$$d\,\E\left(\int_0^t\frac{ds}{\rho_{2,s}}\right)=(d-1)\E\left(p.v.\int_0^t\frac{ds}{\rho_{1,s}}\right).$$
Because of Brownian scaling, the two sides are proportional to
$t\mapsto\sqrt t$. One can check directly that the coefficients agree. Indeed, from the explicit distribution of
$\rho_{1,t}$, $\rho_{2,t}$ (starting from 0), it is easy to compute:
$$\E(\rho_{1,t})=\sqrt{2t}\frac{\Gamma\left(\frac{d+1}2\right)}{\Gamma\left(\frac{d}2\right)},
\ \ \ \ \ \ \ \ 
\E(\rho_{2,t})=\sqrt{2t}\frac{\Gamma\left(\frac{d+3}2\right)}{\Gamma\left(\frac{d+2}2\right)}.
$$
From the Dirichlet decompositions, one gets:
\begin{align*}
(d-1)\E\left(p.v.\int_0^t\frac{ds}{\rho_{1,s}}\right)=2\E(\rho_{1,t})=2\frac{\Gamma\left(\frac{d+1}2\right)}{\Gamma\left(\frac{d}2\right)}\sqrt{2t}\\
d\,\E\left(\int_0^t\frac{ds}{\rho_{2,s}}\right)=\frac{2d}{d+1}\E(\rho_{2,t})=\frac{2d}{d+1}.\frac{\Gamma\left(\frac{d+3}2\right)}{\Gamma\left(\frac{d+2}2\right)}\sqrt{2t}=2\frac{\Gamma\left(\frac{d+1}2\right)}{\Gamma\left(\frac{d}2\right)}\sqrt{2t}.
\end{align*}

In the limiting case $d\nearrow 1$, $\rho_1$ is a reflecting Brownian
motion that can be represented as: $|B_t|=\tilde B_t+\ell^0_t$, where
$\ell^0$ is the local time at 0 of the Brownian motion
$B$ and $\tilde B$ is also a Brownian motion. Assuming that the Dirichlet decomposition of $\rho_1$,
i.e. $(B^1_t,\frac{d-1}2\int_0^tds/\rho_{1,s})$ converges to the
Dirichlet (or semimartingale) decomposition of $|B|$, that is $(\tilde
B_t,\ell^0_t)$, then (iv) translates into:
$$\tilde B_t+2\ell^0_t=\rho_{2,t}$$
where $\rho_2$ is a 3-dimensional Bessel process; this is Pitman's
$(2M-X)$ theorem. One also gets the intertwining relation
$$\E(f(\rho_{1,t})|{\mc F}^2_t)=\int_0^1f(u\rho_{2,t})du$$
which is used in Pitman's original proof.

Let us get back to $\SLE(\kappa,\kappa-4)$. In the discussion
in \cite{D4}, Section 5, we mentioned that the future of an
$\SLE(\kappa,\kappa-4)$ after a random time where the trace lies on
the (final) right boundary is a translated
$\SLE_{(0,0^+)}(\kappa,\kappa-4)$. Since the trace lies on the right
boundary exactly when the process $X$ vanish, a rigorous statement and
proof were given in the proof of Proposition \ref{BP3}. In particular,
we consider the stopping times $\tau_u$, $u\geq 0$, where $\tau$ is
the right-continuous inverse of the local time of $X$ at $0$. These are ${\mc
  X}$-stopping times a.s. supported on the zero set of $X$, so that if
$(K_t)$, $(g_t)$ denote the families of hulls and conformal
equivalences respectively, then for any $u\geq 0$,
$(g_{\tau_u}(K_{\tau_u+t}\setminus K_{\tau_u})-W_{\tau_u})_t$ defines an
$\SLE_{(0,0^+)}(\kappa,\kappa-4)$ independent from $K_{\tau_u}$. 

Let ${\mc Q}$ be the set of bounded hulls in $\H$; recall that a
bounded hull $A$ in $\H$ is a compact subset of $\overline\H$ such
that $\H\backslash A$ is simply connected and $A\cap
\R\subset\overline{A\cap\H}$. Then ${\mc Q}\times \R$ can be seen as a
semigroup of ``pointed hulls'' with composition law:
$$(A,x).(B,y)=(A\cup\phi_A^{-1}(B+x),x+y)$$
Note that $(A,x)\mapsto (\capp(A),x)$, where $\capp$ denotes the
half-plane capacity, is a semigroup morphism ${\mc
  Q}\times\R\rightarrow \R^+\times \R$. The only invertible elements in ${\mc Q}\times \R$ are the
$(\varnothing,x)$, $x\in\R$.

Consider now the ${\mc Q}\times \R$-valued process: $u\mapsto
(K_{\tau_u},W_{\tau_u})$. Then this process has the independent
increment property; this follows readily from the fact that
$(g_{\tau_u}(K_{\tau_u+t}\setminus K_{\tau_u})-W_{\tau_u})_t$ defines an
$\SLE_{(0,0^+)}(\kappa,\kappa-4)$ independent from $K_{\tau_u}$.
 Taking into account the restriction formulae discussed
in the previous chapter, one can formulate the:

\begin{Prop}
Let $H_u=K_{\tau_u}$, $w_u=W_{\tau_u}$. Then $(H_u,w_u)$ is a
${\mc Q}\times\R$-valued stable L\'evy process with index $\nu=\frac 12+\frac
2\kappa$ in the following sense:\\
(i) For any $u,v\geq 0$, if $(\tilde H,\tilde w)$ is an independent
copy of $(H,w)$, then:
$$(H_u,w_u).(\tilde H_v,\tilde w_v)\stackrel{\mc L}=(H_{u+v},w_{u+v}).$$
(ii) For any $c>0$, $u\geq 0$, $(H_{cu},w_{cu})_u\stackrel{\mc
  L}{=}c^{1/2\nu}(H_u,w_u)_u$.\\
Moreover, $(H_u,w_u)$ satisfies the following restriction formula:
for any smooth $+$-hull $A$, if $\alpha=\frac 12-\frac 1\kappa$, and $L$ is
an independent loop soup with intensity $\lambda_\kappa$, then:
$$\phi'_A(0)^\alpha=\E\left(\phi'_{\phi_{H_u}(A)}(w_u)^\alpha\ind_{H_u\cap
    A^L=\varnothing}\right).$$
\end{Prop}

Obviously, the image of a semigroup-valued L\'evy process under a
semigroup homomorphism is again a L\'evy process. Of particular
interest is the $\R^+\times\R$-valued process:
$$(\tau_u,w_u)_u=\left(\tau_u,p.v.\int_0^{\tau_u}\frac{2ds}{X_s}\right)$$
where the second coordinate is a well-known stable L\'evy process with index $2\nu=1+4/\kappa$.

We now describe the excursion process for these excursions ``away from
the right-boundary''. Consider the space $U$ of real continuous
processes with finite lifetime (with the usual topology): an element $w$ of $U$ defines a
continuous real-valued function on $[0,T(w)]$, where $T(w)>0$ is the
lifetime of $w$. Then one can define a Loewner map from a
subset of $U$ to ${\mc Q}\times \R$ in the following way: for a
suitable $w\in U$, with lifetime $T>0$, let $(\tilde g_t)$ be the Loewner flow associated
with the driving function $w$, and $(\tilde K_t)$ the corresponding hulls. Then
${\bf L}w=(\tilde K_T,w_T)$ defines an element of
${\mc Q}\times\R$.

We recall a representation of  excursions of Bessel processes (see
e.g. \cite{PY84}; see also \cite{Ber90}). The excursion measure $n_d$
of a $d$-dimensional Bessel process, $0<d<2$,
seen as a measure on $U$, can be described as follows (up to a
multiplicative constant):
\begin{itemize}
\item The maximum $M$ of the excursion satisfies $n_d(M>x)=x^{d-2}$
  for $x>0$.
\item Under the probability measure $n_d(.\ind_{M>x})/n_d(M>x)$, the
  excursion is the concatenation of a $(4-d)$-dimensional Bessel
  process, run until it hits the level $x$, and an independent
  $\BES(d)$ process started from $x$, killed when it hits 0.
\end{itemize}

We briefly justify this decomposition. If $T_x$ is the first time at
which the $d$-dimensional Bessel process $\rho$ (started from 0) reaches the level $x$,
consider the excursion straddling $T_x$. From the Poisson process
property, this excursion has distribution $n_d(.|M\geq x)$; so the only thing to check is that $(\rho_u)_{g_{T_x}\leq u\leq
  T_x}$ is a $\BES(4-d)$ stopped when it hits $x$. This follows from
the realization of the $\BES(4-d)$ process as a Doob $h$-transform of
the $\BES(d)$ process.

This decomposition translates into a description of the excursions
``away from the boundary'' for $\SLE(\kappa,\kappa-4)$. We retain the
previous notations.

\begin{Prop}
The process $(e_u)_u$ defined by $e_u=(g_{\tau_{u^-}}^{-1}(K_{\tau_u}\setminus
K_{\tau_{u^-}}), W_{\tau_u}- W_{\tau_{u^-}})$ if $\tau_u>\tau_{u^-}$,
and by $e_u=\partial$ if $\tau_u=\tau_{u^-}$, is a $({\mc Q}\times
\R)\cup\{\partial\}$-valued $({\mc X}_{\tau_u})$- Poisson point process. Its characteristic
measure $N_\kappa$ is the image of $n_d$ under the composition:
$$\begin{array}{ccccc}
U&\longrightarrow&U&\longrightarrow&{\mc Q}\times \R\\
(v_t)_{t<T}&\longmapsto& \left(-\sqrt\kappa
v_t+2\int_0^t\frac{ds}{\sqrt\kappa v_s}\right)_{t<T}\\
&&w&\longmapsto& {\bf L}w
\end{array}
$$
(with the natural extension for cemetary states). In other words, the
excursion measure $N_\kappa$ can be described (up to a multiplicative
constant) as the monotone limit of the measures $y^{d-2}N_{\kappa,y}$
as $y\searrow 0$, where $N_{\kappa,y}$ is the probability measure
resulting from the concatenation of:
\begin{itemize}
\item an $\SLE_{(0,0^+)}(\kappa,\kappa)$, with associated driving
  process $\hat W$ and conformal equivalences $\hat g_t$, run until
  $(\hat g_t(0^+)-\hat W_t)$ hits $y$ at time $\hat T$
\item and an independent (conditionally on its initial state)
  $\SLE_{(\hat W_{\hat T},\hat g_{\hat T}(0^+))}(\kappa,-4)$ with associated driving
  process $\tilde W$ and conformal equivalences $\tilde g_t$, run until
  $(\tilde g_t(\hat W_{\hat T}+y)-\tilde W_t)$ hits $0$.
\end{itemize}
\end{Prop}

Recall the notations of \cite{D4}, Section 6. A (somewhat wishful) roadmap to proving duality would involve:
\begin{itemize}
\item Computing the tails of the distribution of the first two moments
  (i.e. translation at infinity and capacity seen from infinity) 
  of the ``synthetic'' hull $K_2$: $\phi_{K_2}(z)=z-w+2\tau/z+O(1/z^2)$, 
\item Proving the convergence in distribution:
 $$\left(\frac{\tau^{(1)}+\cdots+\tau^{(n)}}{n^{1/\nu}},\frac{w^{(1)}+\cdots+w^{(n)}}{n^{1/2\nu}}\right)\longrightarrow (\tau_l,\int_0^{\tau_l}2ds/X_s)$$
where $(\tau^{(i)},w^{(i)})$ are the first two moments of an independent
copy $K_2^{(i)}$ of $K_2$ and $l$ is some constant,
\item Proving that a stable law on hulls is determined by the joint
  law of its first two moments under suitable conditions, and
\item Proving that the right boundary of the concatenation $K_2^{(1)}\dots K_2^{(n)}$ is an
  $\SLE_l(\kappa',\kappa'/2-2)$ stopped at a random time.
\end{itemize}

\section{Cutpoints for $\SLE(\kappa,\kappa-4,\kappa-4)$}

In the previous section, we were able to explicit the different
processes involved, since they were $\R^+$-valued diffusions satisfying
Brownian scaling, hence Bessel processes. In the two-sided case we are
now discussing, one has to recast the arguments in a more abstract
fashion, that also sheds new light on the one-sided case studied in the previous section. We
begin with a rephrasing of the previous construction, based on Doob
$h$-transforms and related Girsanov transformations (see \cite{W2},
where Girsanov densities are interpreted in terms of restriction measures).

As before, $\kappa>4$ is fixed, $d=1-4/\kappa$; $(W_t)$ is the Brownian
motion driving an $\SLE_0(\kappa)$ process, and $g_t$ are the
associated conformal equivalences. The notations $\P$ and $\E$ refer to
this process.

For any $y>0$, $Y_t=g_t(y)-W_t$ defines a $\BES(2-d)$ process; as is
well known, $h(u)=u^d$ is a scale function for this diffusion. Hence,
for any $y>0$ and $t>0$, if $\tau_y$ denotes the swallowing time of
$y$ and $Y_t=g_t(y)-W_t$, then:
$$(P^{\uparrow}_tf)(y)=\frac 1{h(y)}\E\left(f(Y_t)\ind_{\tau_y>t}h(Y_t)\right)$$
defines a Markov semigroup (\,$(\ind_{\tau_y>t}h(Y_t))_t$ is an $L^1$-
martingale under $\P$). Note that the probability measures $P_t(y,.)$
admit a weak limit as $y\searrow 0$. The infinitesimal generator of
this semigroup is obtained by conjugation
of the original generator, and it is easily seen that
$(P^{\uparrow}_t)$ is the semigroup of $(Y_t=g_t(y)-W_t)$ under
the $\SLE(\kappa,\kappa-4)$ probability measure. Let $\P^\uparrow$ be a
corresponding probability measure (i.e. under $\P^\uparrow$,
$Y/\sqrt\kappa$ is a $\BES(2+d)$ process). In terms of Bessel
processes, this is a classical result, see e.g. \cite{PY84}.

Now, let $0<x<y$, $Y_t=g_t(y)-W_t$ as above, and
$X_t=g_t(x)-W_t$. We look at $(X_t)$ as a functional of $(Y_t)$, as a
consequence of the identity:
$$d(X_t-Y_t)=-\frac{2(X_t-Y_t)dt}{X_tY_t}.$$
Then the semimartingale $(\ind_{t<\tau_x}h(X_t)/h(Y_t))$ is a
$\P^\uparrow$-martingale. Indeed
$$\frac 1{h(y)}\E\left(\left(\frac {h(X_t)}{h(Y_t)}\ind_{\tau_x>t}\right)\ind_{\tau_y>t}h(Y_t)\right)=\frac 1{h(y)}\E\left(\ind_{\tau_x>t}h(X_t)\right)=\frac{h(x)}{h(y)}.$$
Furthermore, $h(x)/h(y)=\P^\uparrow(\tau_x=\infty)$, so one gets for
any bounded Borel function $f$:
$$\E^\uparrow_y(f(X_t)|\tau_x=\infty)=\frac{h(y)}{h(x)}.\frac 1{h(y)}\E\left(f(X_t)\left(\frac
    {h(X_t)}{h(Y_t)}\ind_{\tau_x>t}\right)\ind_{\tau_y>t}h(Y_t)\right)=P^\uparrow_t f(x),$$
which is essentially Proposition \ref{BP1} (i)-(ii). This is formally
very similar to the following properties of the three-dimensional
Bessel process: the $\BES(3)$ process is the Doob $h$-transform of a Brownian
motion, with $h(x)=x$; if $\rho$ is a $\BES(3)$ process, then
conditionally on $\inf_{s\geq 0}(\rho_s)\geq a$, $\rho-a$ is again a
$\BES(3)$ process.

We have also used a description of $\SLE(\kappa,\kappa-4)$
conditionally on $Z$, its rightmost swallowed point. In the present set-up,
note that $1-h(x)/h(y)=\P^\uparrow(\tau_x<\infty)$, so:
\begin{align*}
\E^\uparrow_y(f(X_t)\ind_{\tau_x>t}|\tau_x<\infty)&=\left(1-\frac{h(x)}{h(y)}\right)^{-1}\frac 1{h(y)}\E\left(f(X_t)\ind_{\tau_x>t}\left(1-\frac
    {h(X_t)}{h(Y_t)}\right)\ind_{\tau_y>t}h(Y_t)\right)\\
&=\frac 1{h(y)-h(x)}\E\left(f(X_t)\ind_{\tau_x>t}(h(Y_t)-h(X_t))\right).
\end{align*}
Now let $x\nearrow y$ (a rigorous proof of these results has already
been given, so we can afford to argue a bit loosely here). One gets:
$$\E^\uparrow_y(f(Y_t)\ind_{\tau_y>t}|Z\in dy)=\frac 1{h'(y)}\E\left(f(Y_t)\ind
_{\tau_y>t}g'_t(y)h'(Y_t)\right)dy.$$
It is indeed easy to check that
$(g'_t(y)(g_t(y)-W_t)^{d-1}\ind_{\tau_y\geq t})_t$ is an
integrable $\P$-local martingale. So we can define a (non conservative)
Markov semigroup:
$$(P^\downarrow_tf)(y)=\frac 1{h'(y)}\E\left(f(Y_t)\ind
_{\tau_y>t}g'_t(y)h'(Y_t)\right)$$
which corresponds to an $\SLE_{(0,y)}(\kappa,-4)$ killed at $\tau_y$.

Finally, the excursion measure of $X$ can be described in terms of
$\BES(4-d)$ and $\BES(d)$ processes. The $\BES(d)$ process corresponds
to the semigroup $P^{\downarrow}$, while the $\BES(4-d)$ process
corresponds to the semigroup $P^{\updownarrow}$, where:
$$(P^\updownarrow_tf)(y)=\frac 1{y}\E\left(f(Y_t)\ind
_{\tau_y>t}g'_t(y)Y_t\right).$$

To sum up, we have encountered $\SLE(\kappa,\rho)$ processes with
$\rho=0,\kappa-4,-4,\kappa$, associated with Bessel processes of
dimensions $\delta=2-d,2+d,d,4-d$ (or index $\nu=-1/2+2/\kappa,1/2-2/\kappa,
-1/2-2/\kappa,1/2+2/\kappa$); the corresponding semigroups are
$P,P^\uparrow,P^\downarrow,P^\updownarrow$, defined by means of the
$\P$- local martingales $(\ind),(Y_t^d),(g'_t(y)Y_t^{d-1}),(g'_t(y)Y_t)$. Note
that these densities are also considered in \cite{W2}.

We now turn to the two-sided case. Under $\P$, $W/\sqrt\kappa$ is a
standard Brownian motion, $g_t$ are the associated conformal
equivalences, and ${\mc F}$ is the natural filtration. If $y_1<0<y_2$, then it is easy to check that:
$$\left(W_t-g_t(y_1)\right)^d\left(g_t(y_2)-W_t\right)^d\left(g_t(y_2)-g_t(y_1)\right)^{\frac\kappa 2d^2}\ind_{\tau_{y_1}>t,\tau_{y_2}>t}$$
is a martingale. Let $k(u_1,u_2)=(-u_1)^du_2^d(u_2-u_1)^{\kappa d^2/2}$,
and define:
$$(Q_t^{\uparrow\uparrow}f) (y_1,y_2)=\frac 1{k(y_1,y_2)}\E\left((fk)(g_t(y_1)-W_t,g_t(y_2)-W_t)
\ind_{\tau_{y_1}>t,\tau_{y_2}>t}\right)$$
for any bounded Borel function $f$ on $\R^-\times\R^+$. As before, this is the semigroup of an $(\R^+)^2$-valued Markov
process. If ${\mc L}$ and ${\mc L}^{\uparrow\uparrow}$ designate the
infinitesimal generator of the diffusion $(g_t(y_1)-W_t,g_t(y_2)-W_t)$
under $\P$ and the generator of $Q^{\uparrow\uparrow}$ respectively, say
restricted to $C^2_0(\R^-\times\R^+)$, then:
\begin{align*}
{\mc L}&=\frac\kappa 2\left(\frac\partial{\partial
    u_1}+\frac\partial{\partial u_2}\right)^2+\frac 2{u_1}\frac\partial{\partial
    u_1}+\frac 2{u_2}\frac\partial{\partial
    u_1}\\
{\mc L}^{\uparrow\uparrow}&=k(u_1,u_2)^{-1}{\mc L}k(u_1,u_2)\\
&=\frac\kappa 2\left(\frac\partial{\partial
    u_1}+\frac\partial{\partial u_2}\right)^2+\left(\frac {\kappa-2}{u_1}+\frac{\kappa-4}{u_2}\right) \frac\partial{\partial
    u_1}+\left(\frac {\kappa-4}{u_1}+\frac{\kappa-2}{u_2}\right)\frac\partial{\partial
    u_2}
\end{align*}
where $k$ is seen as a multiplication operator. So ${\mc L}^{\uparrow\uparrow}$ is the
generator of an $\SLE(\kappa,\kappa-4,\kappa-4)$, and we have
identified $Q^{\uparrow\uparrow}$. If $y_1<x_1<0<x_2<y_2$, it now appears that
$$\frac{k\left(g_t(x_1)-W_t,g_t(x_2)-W_t\right)}{k\left(g_t(y_1)-W_t,g_t(y_2)-W_t\right)}\ind_{\tau_{x_1}>t,\tau_{x_2}>t}$$
is a martingale under an
$\SLE_{(0,y_1,y_2)}(\kappa,\kappa-4,\kappa-4)$ probability measure,
which is Proposition \ref{BP2} (i). Proposition \ref{BP2} (ii) also follows immediately.

We now discuss the two-sided analogue of Proposition \ref{BP1} (iv),
which is more or less straightforward if a bit tedious. Note that
claims on regular conditional probabilities can be made rigorous,
along the lines of the proof of Proposition \ref{BP1}.
Let $Z_1$ (resp. $Z_2$) the
leftmost (resp. rightmost) point swallowed by an
$\SLE_{(0,y_1,y_2)}(\kappa,\kappa-4,\kappa-4)$. Then $(Z_1,Z_2)$ has
density:
$$-\frac{\partial^2}{\partial u_1\partial
  u_2}k(u_1,u_2)\ind_{y_1\leq u_1\leq 0\leq u_2\leq y_2}du_1du_2/k(y_1,y_2).$$
Conditioning on $\{u_1\leq Z_1,Z_2\leq u_2\}$, we get an
$\SLE_{(0,u_1,u_2)}(\kappa,\kappa-4,\kappa-4)$. As in the one-sided
case, if we define $U^i_t=g_t(u_i)-W_t$, $V^i_t=g_t(v_i)-W_t$,
$i=1,2$, it follows that:
%If $f$ is any bounded
%Borel function on $\R^-\times\R^+$, one gets:
$$\E^{\uparrow\uparrow}\left(f(V^1_t,V^2_t)\ind_{\tau_{v_1}>t,\tau_{v_2}>t}\left|\begin{array}{c}
 u_1\leq Z_1\leq
v_1\\v_2\leq Z_2\leq u_2\end{array}\right.\right)=
\frac{
\E\left(f(V^1_t,V^2_t)\ind_{\tau_{v_1}>t,\tau_{v_2}>t}\tilde
    k(U^1_t,U^2_t,V^1_t,V^2_t)\right)}{\tilde k(u_1,u_2,v_1,v_2)}$$
where $\tilde k$ is defined by:
\begin{align*}
\tilde k(u_1,u_2,v_1,v_2)&\stackrel{\rm
  def}=k(u_1,u_2)-k(v_1,u_2)-k(u_1,v_2)+k(v_1,v_2)\\
&=\partial_{u_1u_2}k(u_1,u_2)(u_1-v_1)(u_2-v_2)+o((u_1-v_1)(u_2-v_2))
\end{align*}
so that $\P_{(u_1,u_2)}^{\uparrow\uparrow}(Z_1\leq v_1,Z_2\geq
v_2)=\tilde k(u_1,u_2,v_1,v_2)/k(u_1,u_2)$.

Taking the limit as $v_1\searrow u_1$, $v_2\nearrow u_2$, one gets a
regular conditional probability for the process stopped at $\tau_{Z_1}\wedge\tau_{Z_2}$:
\begin{align*}
\lefteqn{\E^{\uparrow\uparrow}(f(U^1_t,U^2_t)\ind_{\tau_{u_1}>t,\tau_{u_2}>t}|Z_1\in
  du_1,Z_2\in du_2)=}\\
&\hspace{1in}\frac{\E\left(f(U^1_t,U^2_t)\ind_{\tau_{Z_1>t},\tau_{Z_2>t}}(-g'_t(u_1)g'_t(u_2)\partial_{u_1u_2}k(U^1_t,U^2_t))\right)}{-\partial_{u_1u_2}k(u_1,u_2)}du_1du_2.
\end{align*}
Let $Q^{\downarrow\downarrow}$ be the semigroup of this conditional
process killed at $\tau_{Z_1}\wedge\tau_{Z_2}$. Suppose that $\tau_{Z_1}<\tau_{Z_2}$; between $\tau_{Z_1}$ and
$\tau_{Z_2}$, the law of the process
$(Z^1_t,Z^2_t)=(g_t(Z_1^-)-W_t,g_t(Z_2^+)-W_t)$ is given by:
\begin{align*}
\lefteqn{\E^{\uparrow\uparrow}\left(f(Z^1_{\tau_{Z_1}+t},Z^2_{\tau_{Z_1}+t})\ind_{\tau_{Z_2}>\tau_{Z_1}+t}|{\mc
  F}_{\tau_{Z_1}},Z_1\in dv_1,Z_2\in
dv_2\right)
=}\\
&\hspace{2in}\frac{\E\left(f(U^1_t,U^2_t)\ind_{\tau_{u_2}>t}g'_t(u_2)\partial_{u_2}k(U^1_t,U^2_t)\right)}{\partial_{u_2}k(u_1,u_2)}dv_1dv_2
\end{align*}
where $u_i=Z^i_{\tau_{Z_1}\wedge \tau_{Z_2}}$, $U^i_t=g_t(u^i)-W_t$,
$i=1,2$. Symmetrically, if $\tau_{Z_2}<\tau_{Z_1}$, then the law of the process between $\tau_{Z_2}$ and
$\tau_{Z_1}$ is given by:
\begin{align*}
\lefteqn{\E^{\uparrow\uparrow}\left(f(U^1_{\tau_{Z_2}+t},U^2_{\tau_{Z_2}+t})\ind_{\tau_{Z_1}>\tau_{Z_2}+t}|{\mc
  F}_{\tau_{Z_2}},Z_1\in dv_1,Z_2\in
dv_2\right)=}\\
&\hspace{2in}\frac{\E\left(f(U^1_t,U^2_t)\ind_{\tau_{u_1}>t}g'_t(u_1)\partial_{u_1}k(U^1_t,U^2_t)\right)}{\partial_{u_1}k(u_1,u_2)}dv_1dv_2.
\end{align*}

We sum up this description of an
$\SLE_{(0,y_1,y_2)}(\kappa,\kappa-4,\kappa-4)$ conditionally on
$(Z_1,Z_2)=(u_1,u_2)$ (in the regular conditional probability
sense). Specifically, we give a path decomposition for $(Z^1,Z^2)$,
where $Z_t^i=g_t(Z_i)-W_t$; $f$ is some bounded Borel function on
$\R^-\times\R^+$, and $\tau_i=\tau_{Z_i}$.
\begin{itemize}
\item Up to time $\tau_1\wedge\tau_2$, $(Z^1_s,Z^2_s)_{s\leq \tau_1\wedge\tau_2}$ is a Markov
  process starting from $(Z_1,Z_2)$ with semigroup:
$$Q_t^{\downarrow\downarrow}f(u_1,u_2)=\frac{\E\left(f(U^1_t,U^2_t)\ind_{\tau_{u_1}>t,\tau_{u_2}>t}(-g'_t(u_1)g'_t(u_2)\partial_{u_1u_2}k(U^1_t,U^2_t))\right)}{-\partial_{u_1u_2}k(u_1,u_2)}.$$
\item If $\tau_1<\tau_2$, then
  $(Z^1_{s+\tau_1},Z^2_{s+\tau_1})_{s\leq\tau_2-\tau_1}$ is a
  Markov process starting from $(0^-,Z^2_{\tau_1})$ with semigroup:
$$Q_t^{\uparrow\downarrow}f(u_1,u_2)=\frac{\E\left(f(U^1_t, U^2_t)\ind_{\tau_{u_1}>t,\tau_{u_2}>t}(g'_t(u_2)\partial_{u_2}k(U^1_t,U^2_t))\right)}{\partial_{u_2}k(u_1,u_2)}.$$
\item If $\tau_2<\tau_1$, then
  $(Z^1_{s+\tau_2},Z^2_{s+\tau_2})_{s\leq\tau_1-\tau_2}$ is a
  Markov process starting from $(Z^1_{\tau_2},0^+)$ with semigroup:
$$Q_t^{\downarrow\uparrow}f(u_1,u_2)=\frac{\E\left(f(U^1_t, U^2_t)\ind_{\tau_{u_1}>t,\tau_{u_2}>t}(-g'_t(u_1)\partial_{u_1}k(U^1_t,U^2_t))\right)}{-\partial_{u_1}k(u_1,u_2)}.$$
\item After time $\tau=\tau_1\vee\tau_2$, the process
  $(Z^1_{\tau+s},Z^2_{\tau+s})_s$ is a Markov process starting from
  $(Z^1_\tau,Z^2_\tau)$ with semigroup:
$$Q_t^{\uparrow\uparrow}f(u_1,u_2)=\frac{\E\left(f(U^1_t,U^2_t)\ind_{\tau_{u_1}>t,\tau_{u_2}>t}(k(U^1_t,U^2_t))\right)}{k(u_1,u_2)}.$$
\end{itemize}

We now consider the following situation: $(g_t)$ (resp. $(K_t)$) is the family of
conformal equivalences (resp. hulls) associated with an
$\SLE_{(0,y_1,y_2)}(\kappa,\kappa-4,\kappa-4)$, $y_1\leq 0\leq y_2$, $W$ is the driving
process of this $\SLE$, $Y^i_t=g_t(y_i)-W_t$. For any time $t$,
let $X^1_t$ (resp. $X^2_t$) be the leftmost (resp. rightmost) point
swallowed by $(g_t^{-1}(K_{t+s}\setminus K_t)-W_t)_s$, which is an
$\SLE_{(0,Y^1_t,Y^2_t)}(\kappa,\kappa-4,\kappa-4)$ process. Alternatively, one
  may define:
$$(X^1_t,X^2_t)=(g_t(\gamma_{\sigma^1_t})-W_t,g_t(\gamma_{\sigma^2_t})-W_t)$$
where $\sigma^1_t$ (resp. $\sigma^2_t$) is the last time before $t$
spent by the trace of the left (resp. right) boundary of the final
hull $K_\infty$. Cuttimes can be characterized as zeroes of the joint
process $(X^1,X^2)$; indeed, if $X^1$ and $X^2$ vanish at time $t$,
then by definition $\gamma_t$ lies on the left and right boundaries of
$K_\infty$, in other words $\gamma_t$ is a cutpoint of $K_\infty$.

We have seen that $(X^1,X^2)$ is a continuous process, and that for
$t>0$, the law of $(X_t^1,X_t^2)$ conditionally on ${\mc
  Y}_t=\sigma((Y^1_s,Y^2_s)_{0\leq s\leq t})$ is
$\Upsilon((Y^1_t,Y^2_t),.)$, where $\Upsilon$ is the transition kernel
$\R^-\times\R^+\rightarrow\R^-\times\R^+$ specified by:
$$(\Upsilon f)(y_1,y_2)=\int_{y_1}^0\int_0^{y_2}f(x_1,x_2)\frac
{-\partial_{x_1x_2}k(x_1,x_2)}{k(y_1,y_2)}dx_1dx_2$$
where $f$ is some bounded Borel function on $\R^-\times\R^+$. We will
also need one-sided versions of this intertwining kernel. More
precisely, the conditional law of $(X_t^1,X_t^2)$ given
$(X_t^1,Y_t^2)$ is $\Upsilon^1((X^1_t,Y^2_t),.)$, where:
$$(\Upsilon^1 f)(x_1,y_2)=\int_0^{y_2}f(x_1,x_2)\frac
{-\partial_{x_1x_2}k(x_1,x_2)}{-\partial_{x_1}k(x_1,y_2)}dx_2.$$
Symmetrically, the law of $(X_t^1,X_t^2)$ given
$(Y_t^1,X_t^2)$ is $\Upsilon^2((Y^1_t,X^2_t),.)$, where:
$$(\Upsilon^2 f)(y_1,x_2)=\int_{y_1}^0f(x_1,x_2)\frac
{-\partial_{x_1x_2}k(x_1,x_2)}{\partial_{x_2}k(y_1,x_2)}dx_1.$$

We will prove that $(X^1,X^2)$ is a Markov process; we now describe its
semigroup. We have seen that the law of $(X^1,X^2)$ starting from
$x_1<0<x_2$ and killed at $\tau_{X_1}\wedge\tau_{X_2}$ is that of a
Markov process with (non conservative) semigroup
$Q^{\downarrow\downarrow}$. So the problem is to extend trajectories
of $(X^1,X^2)$ after $\tau_{X_1}\wedge\tau_{X_2}$; for this, we
will need the semigroups $Q^{\uparrow\downarrow}$,
$Q^{\downarrow\uparrow}$, $Q^{\uparrow\uparrow}$, and the intertwining
kernels $\Upsilon^2$, $\Upsilon^1$, $\Upsilon$. 

Let $(X^{\downarrow\downarrow}_1,X^{\downarrow\downarrow}_2)$ be a Markov process with
(non conservative) semigroup $Q^{\downarrow\downarrow}$, started from
$(x_1,x_2)$, $x_1\leq 0\leq x_2$, and $T$ be its lifetime; then either
$X^{\downarrow\downarrow}_{1,T}$ or $X^{\downarrow\downarrow}_{2,T}$
vanish. If $X^{\downarrow\downarrow}_{1,T}=0$ (resp. $X^{\downarrow\downarrow}_{2,T}=0$), consider an
independent Markov process
$(Y^{\uparrow\downarrow}_1,X^{\uparrow\downarrow}_2)$
(resp. $(X^{\downarrow\uparrow}_1,Y^{\downarrow\uparrow}_2)$\,) with
semigroup $Q^{\uparrow\downarrow}$ (resp. $Q^{\downarrow\uparrow}$)
started from $(0^-,X^{\downarrow\downarrow}_{2,T})$
(resp. $(X^{\downarrow\downarrow}_{1,T},0^+)$\,); let $T'$ be its
lifetime. Finally, let
$(Y^{\uparrow\uparrow}_1,Y^{\uparrow\uparrow}_2)$ be an independent
Markov process with semigroup $Q^{\uparrow\uparrow}$ started from
$(Y^{\uparrow\downarrow}_1,X^{\uparrow\downarrow}_2)_{T'}$ (resp. 
$(X^{\downarrow\uparrow}_1,Y^{\downarrow\uparrow}_2)_{T'}$\,). These
processes are defined on some probability space; the notations ${\bf
  P}$, ${\bf E}$ refer to this space. One can now define:
\begin{align*}
(\tilde
Q_t^{\downarrow\downarrow}f)(x_1,x_2)=&\hphantom{+}{\bf E}\left(\ind_{t\leq T}f\left(X^{\downarrow\downarrow}_{1,t},X^{\downarrow\downarrow}_{2,t}\right)\right)\\
&+{\bf E}\left(\ind_{T<t\leq
    T+T'}\ind_{X^{\downarrow\downarrow}_{1,T}=0}(\Upsilon^2
  f)\left(Y^{\uparrow\downarrow}_{1,t-T},X^{\uparrow\downarrow}_{2,t-T}\right)\right)\\
&+{\bf E}\left(\ind_{T<t\leq
    T+T'}\ind_{X^{\downarrow\downarrow}_{2,T}=0}(\Upsilon^1
  f)\left(X^{\downarrow\uparrow}_{1,t-T},Y^{\downarrow\uparrow}_{2,t-T}\right)\right)\\
&+{\bf E}\left(\ind_{T+T'<t}(\Upsilon
  f)\left(Y^{\uparrow\uparrow}_{1,t-T-T'},Y^{\uparrow\uparrow}_{2,t-T-T'}\right)\right).
\end{align*}

It will also be convenient to extend the non conservative
semigroups $Q^{\uparrow\downarrow}$, $Q^{\downarrow\uparrow}$. Denote
by $\iota$ the following involution of $\R^-\times \R^+$:
$\iota((u_1,u_2))=(-u_2,-u_1)$. Then $Q^{\uparrow\downarrow}=\iota
Q^{\downarrow\uparrow}\iota$, $\Upsilon^1=\iota \Upsilon^2\iota$, so
that we can restrict to $Q^{\uparrow\downarrow}$ for instance. As
above, we consider a Markov process
$(Y^{\uparrow\downarrow}_1,X^{\uparrow\downarrow}_2)$ with semigroup
$Q^{\uparrow\downarrow}$ started from $(y_1\leq 0\leq x_2)$; let $T$
be its lifetime (i.e. $X^{\uparrow\downarrow}_2$ is positive before
$T$ and vanishes at time $T$). Let
$(Y^{\uparrow\uparrow}_1,Y^{\uparrow\uparrow}_2)$ be a Markov process
with semigroup $Q^{\uparrow\uparrow}$ independent of the former
conditionally on its starting state
$(Y^{\uparrow\downarrow}_{1,T},0^+)$. We will also need another two Markov
kernels (from $\R^-\times\R^+$ to itself):
\begin{align*}
(\Upsilon_1f)(y_1,y_2)&=\int_0^{y_2}f(y_1,x_2)\frac{\partial_{x_2}k(y_1,x_2)}{k(y_1,y_2)}dx_2\\
(\Upsilon_2f)(y_1,y_2)&=\int_{y_1}^0f(x_1,y_2)\frac{-\partial_{x_1}k(x_1,y_2)}{k(y_1,y_2)}dx_1
\end{align*}
so that the conditional law of $(Y^1_t,X^2_t)$ given $(Y^1_t,Y^2_t)$
is $\Upsilon_1((Y^1_t,Y^2_t),.)$, and the conditional law of $(X^1_t,Y^2_t)$ given $(Y^1_t,Y^2_t)$
is $\Upsilon_2((Y^1_t,Y^2_t),.)$; the identities
$\Upsilon_1\Upsilon^2=\Upsilon_2\Upsilon^1=\Upsilon$ are then
obvious. We now define:
\begin{align*}
(\tilde
Q_t^{\uparrow\downarrow}f)(y_1,x_2)=&{\bf E}\left(\ind_{t\leq
    T}f\left(Y^{\uparrow\downarrow}_{1,t},X^{\uparrow\downarrow}_{2,t}\right)\right) 
+{\bf E}\left(\ind_{T<t}(\Upsilon_1 
  f)\left(Y^{\uparrow\uparrow}_{1,t-T},Y^{\uparrow\uparrow}_{2,t-T}\right)\right).
\end{align*}

Note that the semigroup properties $\tilde Q_{t+s}^{\downarrow\downarrow}=\tilde
Q_t^{\downarrow\downarrow}\tilde
Q_s^{\downarrow\downarrow}$, $\tilde Q_{t+s}^{\uparrow\downarrow}=\tilde
Q_t^{\uparrow\downarrow}\tilde
Q_s^{\uparrow\downarrow}$, are not obvious for the moment. Let ${\mc
  X}^i_t=\sigma((X^i_s)_{0\leq s\leq t})$, ${\mc X}={\mc X}^1\vee{\mc X}^2$; we can finally state:

\begin{Prop} The process $(X^1,X^2)$ (resp. $(Y^1,X^2)$, $(X^1,Y^2)$) is Markov with semigroup $\tilde
  Q^{\downarrow\downarrow}$ (resp. $\tilde
  Q^{\uparrow\downarrow}$, $\tilde Q^{\downarrow\uparrow}$) in the filtration ${\mc X}\vee{\mc
    Y}$ (resp. ${\mc X}^2\vee{\mc
    Y}$, ${\mc X}^1\vee{\mc
    Y}$) . Moreover, the following intertwining relations hold:
\begin{align*}
&&Q^{\uparrow\uparrow}\Upsilon&=\Upsilon \tilde Q^{\downarrow\downarrow}\\
Q^{\uparrow\uparrow}\Upsilon_1&=\Upsilon_1 \tilde
Q^{\uparrow\downarrow}&&&
\tilde Q^{\uparrow\downarrow}\Upsilon^2&=\Upsilon^2 \tilde
Q^{\downarrow\downarrow}
\\
Q^{\uparrow\uparrow}\Upsilon_2&=\Upsilon_2 \tilde
Q^{\downarrow\uparrow}&&&\tilde Q^{\downarrow\uparrow}\Upsilon^1&=\Upsilon^1 \tilde
Q^{\downarrow\downarrow}
\end{align*}
 
\end{Prop}
\begin{proof}
First, we consider the statement for $(Y^1,X^2)$ (the case of
$(X^1,Y^2)$ being symmetric). Let $0\leq s\leq t$. As in the one-sided case, we see that $(Y^1_{s+u},X^2_{s+u})_u$ is
independent from ${\mc X}^2_s\vee{\mc Y}_s$ conditionally on
$(Y_s^1,X^2_s)$; this follows from the Markov property of
$\SLE(\kappa,\kappa-4,\kappa-4)$, Proposition \ref{BP2} (ii), and the
fact that ${\mc X}^i_s\vee{\mc Y}_s=\sigma({\mc
    Y}_s,X_s^i)$, $i=1,2$. Moreover, $(Y^1_{s+u},X^2_{s+u})_u$ killed at
  $\tau_2$ is Markov with semigroup
  $Q^{\uparrow\downarrow}$, where $\tau_2$ is the first time after
  $s$ at which $X^2$ vanishes.

Let $\tau$ be an ${\mc X}^2\vee{\mc Y}$ stopping time such that
$X^2_\tau=0$ a.s.\,. Then
$$(\tilde Y^1,\tilde
Y^2)_u=(Y^1_{\tau+u},g_{\tau+u}(g_\tau^{-1}(W_\tau^+))
-W_{\tau+u})_u$$
is a Markov process with semigroup $Q^{\uparrow\uparrow}$,
independent from $({\mc X}^2\vee{\mc Y})_\tau$ conditionally on its
starting state $(Y^1_\tau,0^+)$. As in Proposition \ref{BP3}, this
follows from the decomposition of an $\SLE(\kappa,\kappa-4,\kappa-4)$
conditionally on its rightmost swallowed point. The process
$(\tilde Y^1,\tilde Y^2)_u$ is but
the driving mechanism of an $\SLE(\kappa,\kappa-4,\kappa-4)$
independent from ${\mc X}^2_t\vee{\mc Y}_t$ conditionally on $Y^1_t$; let $(\tilde X^1,\tilde X^2)$ be the associated process
(i.e. $\tilde X^2_u$ is the rightmost swallowed point after time $u$). Then
$\tilde X_u^2=X^2_{\tau+u}$ and the conditional law of $(\tilde
Y^1_u,\tilde X^2_u)$ given $(\tilde
Y^1_u,\tilde Y^2_u)$ is $\Upsilon_1((\tilde
Y^1_u,\tilde Y^2_u),.)$. Setting $\tau=\tau_2$, one gets:
\begin{align*}
\E(f(Y^1_t,X^2_t)|{\mc X}^2_s,{\mc Y}_s)&=
(Q^{\uparrow\downarrow}_{t-s}f)(Y^1_s,X^2_s)+
\E\left(\ind_{\tau_2\leq t}(Q^{\uparrow\uparrow}_{t-\tau_2}\Upsilon_1f)(Y^1_{\tau_2},0^+)\right)\\
&=(\tilde Q_{t-s}^{\uparrow\downarrow} f)(Y^1_s,X^2_s).
\end{align*}
Hence it appears that $(Y^1,X^2)$ is an homogeneous ${\mc X}^2\vee{\mc Y}$-Markov process
with semigroup $\tilde Q_{t-s}^{\uparrow\downarrow}$. As for the
intertwining relation, we get as before:
\begin{align*}
\E(f(Y^1_t,X^2_t)|{\mc Y}_s)&=\E(f(Y^1_t,X^2_t)|{\mc Y}_t|{\mc
  Y}_s)=\E((\Upsilon_1f)(Y^1_t,Y^2_t)|{\mc
  Y}_s)=(Q_{t-s}^{\uparrow\uparrow}\Upsilon_1f)(Y_s^1,Y_s^2)\\
&=\E(f(Y^1_t,X^2_t)|{\mc Y}_s,{\mc X}^2_s|{\mc
  Y}_s)=\E((\tilde Q_{t-s}^{\uparrow\downarrow}f)(Y^1_s,X^2_s)|{\mc
  Y}_s)=(\Upsilon_1\tilde Q_{t-s}^{\uparrow\downarrow}f)(Y_s^1,Y_s^2).
\end{align*}

We now turn to $(X^1,X^2)$, building on the previous result. Let
$0\leq s\leq t$; we have seen that $(X^1_{s+u},X^2_{s+u})_u$ is
independent from ${\mc X}_s\vee{\mc Y}_s$ conditionally on
$(X_s^1,X^2_s)$ (recall that
${\mc X}_s\vee{\mc Y}_s=\sigma({\mc
    Y}_s,X_s^1,X_s^2)$). Moreover, $(X^1_{s+u},X^2_{s+u})_u$ killed at
  $\tau_1\wedge\tau_2$ is Markov with semigroup
  $Q^{\downarrow\downarrow}$, where $\tau_i$ is the first time after
  $s$ at which $X^i$ vanishes.

Suppose that $\tau_1<\tau_2$; then
$(\tilde Y^1_u,\tilde X^2_u)_u=(g_{\tau_1+u}(g_{\tau_1}^{-1}(W_{\tau_1+u}^-))-W_{\tau_1+u},X^2_{\tau_1+u})_u$, killed at
$\tau_2-\tau_1$, is a Markov process with semigroup $\tilde
Q^{\uparrow\downarrow}$, independent from $({\mc X}\vee{\mc Y})_{\tau_1}$
conditionally on its starting state $(0^-,X^2_{\tau_1})$. Moreover,
the law of $(X^1_{\tau_1+u},X^2_{\tau_1+u})$ given $(\tilde
Y^1_u,\tilde X^2_u)$ is $\Upsilon^2((\tilde
Y^1_u,\tilde X^2_u),.)$. It follows that:
\begin{align*}
\E(f(X^1_t,X^2_t)|{\mc X}_s,{\mc
  Y}_s)&=(Q^{\downarrow\downarrow}f)(X^1_s,X^2_s)+\E\left(\ind_{\tau_1\leq
  t\wedge\tau_2}(\tilde Q_{t-\tau_1}^{\uparrow\downarrow}\Upsilon^2f)(0^-,X^2_{\tau_1})\right)\\
&\hphantom{0.3in}+\E\left(\ind_{\tau_2\leq
  t\wedge\tau_1}(\tilde
Q_{t-\tau_2}^{\downarrow\uparrow}\Upsilon^1f)(X^1_{\tau_2},0^+)\right)\\
&=(\tilde Q_{t-s}^{\downarrow\downarrow}f)(X^1_s,X^2_s),
\end{align*}
taking into account the previous decomposition of $\tilde
Q^{\uparrow\downarrow}$, $\tilde
Q^{\downarrow\uparrow}$. The intertwining relations follow as
before. Note that $Q^{\uparrow\uparrow}\Upsilon_1=\Upsilon_1 \tilde
Q^{\uparrow\downarrow}$, $\tilde
Q^{\uparrow\downarrow}\Upsilon^2=\Upsilon^2 \tilde
Q^{\downarrow\downarrow}$, and $\Upsilon=\Upsilon_1\Upsilon^2$ imply
that $Q^{\uparrow\uparrow}\Upsilon=\Upsilon \tilde Q^{\downarrow\downarrow}$. 
\end{proof}

From here, as we can no longer rely on classical results for Bessel
processes, we have to prove the existence of a local time at $(0,0)$
for the Markov process $(X^1,X^2)$, whose right-continuous inverse is
a stable process. The key tool is a local martingale, acting
in many respects as a ``scale function'' for $(X^1,X^2)$.

\begin{Lem}
Let $4<\kappa<8$. If $(X_1,X_2)$ is a Markov process with semigroup
$\tilde Q^{\downarrow\downarrow}$ started from $(x_1,x_2)$, $x_1\leq
0\leq x_2$, then the process:
$$R_t\stackrel{\mathrm {def}}=\frac{(X_{2,t}-X_{1,t})^{4-\kappa/2}}{1+(4-\kappa/2)(X_{1,t}X_{2,t})/(X_{2,t}-X_{1,t})^2}\ind_{t\leq
  T},$$
where $T=\inf(t>0:X_{1,t}=X_{2,t}=0)$, is an ${\mc X}\vee{\mc Y}$- local martingale.
\end{Lem}
\begin{proof}
Note that it is not obvious that $R$ is a semimartingale. It is easily
seen that under $\P$, if $u_1<0<u_2$, the following semimartingale:
$$A_t(u_1,u_2)\stackrel{\mathrm {def}}=g'_t(u_1)g'_t(u_2)(W_t-g_t(u_1))^{d-1}(g_t(u_2)-W_t)^{d-1}(g_t(u_2)-g_t(u_1))^{(\kappa/2)(d-1)^2}$$
is a local martingale (as long as $u_1$ and $u_2$ are not swallowed).
Suppose that $x_1<0<x_2$. If $\tau_i=\inf(t>0:X_{i,t}=0)$, then
$(X_1,X_2)$ killed at $\tau=\tau_1\wedge\tau_2$ is Markov with semigroup
$Q^{\downarrow\downarrow}$. Let $M_n=\inf(t>0:(-X_{1,t})\vee X_{2,t}\geq n)$ be a
  sequence of stopping times. Then:
\begin{align*}
\E^{\downarrow\downarrow}(R_{t\wedge M_n\wedge \tau}|{\mc X_0}\vee{\mc Y}_0)/R_0&=\frac{\E\left(R_ {t\wedge M_n\wedge \tau}(g'_t(x_1)g'_t(x_2)\partial_{u_1u_2}k(X_{1,t},X_{2,t}))\right)}{R_0\partial_{u_1u_2}k(x_1,x_2)}\\
&=\frac{\E(A_{t\wedge M_n\wedge \tau}(x_1,x_2))}{A_0(x_1,x_2)}=1
\end{align*}
which is enough to prove that $(R^\tau_t)$ is a local martingale,
given the Markov property of $X$. Note that we have used the identity:
$$\partial_{u_1u_2}k(u_1,u_2)=d^2(\kappa/2-1)\left(\frac{1}{u_1u_2}+\frac{4-\kappa/2}{(u_2-u_1)^2}\right)k(u_1,u_2).$$ 
We now consider the case $x_1=0<x_2$ (the case $x_1<0=x_2$ being
symmetrical). Suppose that $y_1<0<x_2$. Then $(Y_1,X_2)$ killed at $\tau_2$ is
Markov with semigroup $Q^{\uparrow\downarrow}$, and is intertwined
with $(X_1,X_2)$ via the Markov kernel $\Upsilon^2$. Let $\tilde
M_n=\inf(t>0: (-Y_{1,t})\vee X_{2,t}\geq n)$. From the definitions of
$Q^{\uparrow\downarrow}$ and $\Upsilon^2$, we get:
\begin{align*}
\E^{\downarrow\downarrow}(R_{t\wedge \tilde M_n\wedge \tau_2}|{\mc X}^2_0\vee{\mc Y}_0)&
=\frac{\E\left((\Upsilon^2R)_{t\wedge \tilde M_n\wedge \tau_2}(g'_t(x_2)\partial_{u_2}k(Y_{1,t},X_{2,t}))\right)}{\partial_{u_2}k(y_1,x_2)}
\end{align*}
where we define:
$$(\Upsilon^2R)_t=\int_{Y_{1,t}}^0 r(x_1,X_{2,t})\frac{-\partial_{u_1u_2}k(x_1,X_{2,t})}{\partial_{u_2}k(Y_{1,t},X_{2,t})}dx_1. 
$$
with
$r(x_1,x_2)=(x_2-x_1)^{4-\kappa/2}/(1+(4-\kappa/2)x_1x_2/(x_2-x_1)^2)$.
It follows that:
\begin{align*}
(\Upsilon^2R)_t(g'_t(x_2)\partial_{u_2}k(Y_{1,t},X_{2,t}))
&=\int_{Y_{1,t}}^0
r(x_1,X_{2,t})(-\partial_{u_1u_2}k(x_1,X_{2,t}))g'_t(x_2)dx_1.
%&=d^2(\kappa/2-1)\int_{y_1}^0 A_t(\tilde x_1,x_2)\ind_{t<\tau_{\tilde x_1}}d\tilde x_1
\end{align*}
This is a local martingale; with the change of variable $x_1=g_t(\tilde
x_1)-W_t$, one can see it as an integrated version of
$A_t(.,x_2)$. More precisely, one has to check that:
$$\left({\mc L}-\frac 2{y_2^2}\right)\left(\int_{y_1}^0(-x_1y_2)^{-4/\kappa}(y_2-x_1)^{8/\kappa}dx_1\right)=0$$
where ${\mc L}$ is the differential operator
$(\kappa/2(\partial_{y_1}+\partial_{y_2})^2+2/y_1\partial_{y_1}+2/y_2\partial_{y_2})$.
A simple way to see this is to observe that
\begin{align*}
\partial_{y_1}\left({\mc L}-\frac
  2{y_2^2}\right)\left(\vphantom{\int_{y_1}^0(-x_1y_2)^{-4/\kappa}(y_2-x_1)^{8/\kappa}dx_1}\dots\right)&=
\left({\mc L}-\frac 2{y_1^2}-\frac
  2{y_2^2}\right)\partial_{y_1}\left(\int_{y_1}^0(-x_1y_2)^{-4/\kappa}(y_2-x_1)^{8/\kappa}dx_1\right)\\
&=-\left({\mc L}-\frac 2{y_1^2}-\frac
  2{y_2^2}\right)\left((-y_1y_2)^{-4/\kappa}(y_2-y_1)^{8/\kappa}\right)=0
\end{align*}
and, as $y_1\nearrow 0$, 
\begin{align*}
\left({\mc L}-\frac
  2{y_2^2}\right)\left(\vphantom{\int_{y_1}^0(-x_1y_2)^{-4/\kappa}(y_2-x_1)^{8/\kappa}dx_1}\dots\right)&=
-\left(\frac\kappa 2\partial_{y_1}+\kappa\partial_{y_2}+\frac
  2{y_1}\right)\left((-y_1y_2)^{-4/\kappa}(y_2-y_1)^{8/\kappa}\right)+o(y_1)\\
&=-(-y_1y_2)^{-4/\kappa}(y_2-y_1)^{8/\kappa}\left(-\frac 2{y_1}-\frac
  4{y_2}+\frac 2{y_1}+\frac 4{y_2-y_1}\right)+o(y_1)=o(y_1).
\end{align*}
Given that the process
$((\Upsilon^2R)_t(g'_t(x_2)\partial_{u_2}k(Y_{1,t},X_{2,t})))$ is a
local martingale, we get: 
$$\E^{\downarrow\downarrow}(R_{t\wedge \tilde M_n\wedge
    \tau_2}|{\mc X}^2_0\vee{\mc Y}_0)=\E^{\downarrow\downarrow}(R_0|{\mc X}^2_0\vee{\mc Y}_0)=(\Upsilon^2R)_0(y_1,x_2).$$ 
Letting $y_1\nearrow 0$, this settles the case $x_1=0<x_2$. From the Markov property of $(X^1,X^2)$, we can
  conclude that $(R_t)$ is a local martingale.

Finally, we briefly justify the fact that $A_t(u_1,u_2)$ is
a.s. continuous. Suppose e.g. that $\tau_{u_2}<\tau_{u_1}$; we need
to prove that $A_t$ goes to $0$ as $t\nearrow\tau_{u_2}$. Since
$4<\kappa<8$, there a.s. exists $v_2>u_2$ with
$\tau_{u_2}=\tau_{v_2}$. We know that $g'_t(u_2)(v_2-u_2)\leq
g_t(v_2)-g_t(u_2)$. Moreover, from harmonic measure considerations, it
appears that $(g_t(v_2)-W_t)$ goes to $0$, while the ratio
$(g_t(v_2)-W_t)/(g_t(u_2)-W_t)$ stays bounded, which is enough to conclude.
\end{proof}

Before proceeding with the study of ``excursions away from
cutpoints'', we discuss the limiting case $\kappa=8$ (note that for
$\SLE(8)$, the set of cut-times has a.s. zero Hausdorff dimension, see
\cite{Bef}).

\begin{Lem} For $\kappa=8$, the process:
$$R_t\stackrel{\mathrm {def}}=\left(\log(X_{2,t}-X_{1,t})-\frac{X_{1,t}X_{2,t}}{(X_{2,t}-X_{1,t})^2}\right)\ind_{t\leq
  T},$$
where $T=\inf(t>0:X_{1,t}=X_{2,t}=0)$, is an ${\mc X}\vee{\mc Y}$-
local martingale. As a consequence, $(0,0)$ is a polar point for the
Markov process $(X_1,X_2)$, and the final hull of $\SLE(8,4,4)$ has
a.s. no cutpoints.
\end{Lem}
\begin{proof}
In the previous lemma, we considered the local martingales (for $4<\kappa<8$):
\begin{align*}
R^{(\kappa)}_t&=\frac{(X_{2,t}-X_{1,t})^{4-\kappa/2}}{1+(4-\kappa/2)(X_{1,t}X_{2,t})/(X_{2,t}-X_{1,t})^2}\ind_{t\leq
  T}\\
&=\left(1+(4-\kappa/2)\left(\log(X_{2,t}-X_{1,t})-\frac{X_{1,t}X_{2,t}}{(X_{2,t}-X_{1,t})^2}\right)\right)\ind_{t\leq
  T}+o(\kappa-8),
\end{align*}
which justifies the definition of $R$ here. It can be checked directly
(or using a limiting argument) that for $u_1<0<u_2$, the
semimartingales:
$$A_t(u_1,u_2)=g'_t(u_1)g'_t(u_2)\left(\log(U_{1,t}-U_{2,t})-\frac{U_{1,t}U_{2,t}}{(U_{2,t}-U_{1,t})^2}\right)\frac{U_{2,t}-U_{1,t}}{\sqrt{-U_{1,t}U_{2,t}}}$$
where $U_{i,t}=g_t(u_i)-W_t$, are local martingales under $\P$ (as
long as $u_1,u_2$ are not swallowed). Arguing as in the previous lemma,
we get that $R$ is a local martingale. The fact that $(0,0)$ is polar
for $(X_1,X_2)$ follows immediately. Since cut-times of $\SLE(8,4,4)$
correspond by construction to the zero set of $(X_1,X_2)$, one gets
that $\SLE(8,4,4)$ as a.s. no cutpoint.
\end{proof}

It is easily seen that if $\SLE(8)$ stopped at a finite time had
cutpoints with positive probability, then the final hull of
$\SLE(8,4,4)$ would also have cutpoints with positive
probability. Consequently, $\SLE(8)$ stopped at a finite time has
a.s. no cutpoints. One can prove this directly, building on the
reversibility of the $\SLE(8)$ trace (\cite{Bef}).

We turn back to the case $4<\kappa<8$. Using the above homogeneous ``scale function'', one can define a local
time at $(0,0)$ for $(X_1,X_2)$, and compute the index of its (stable)
inverse, using L\'evy's upcrossings construction. Specifically,
consider the level lines:
$$M_h=\{(x_1,x_2):x_1\leq 0\leq x_2, r(x_1,x_2)=h\}$$
where
$r(x_1,x_2)=(x_2-x_1)^{4-\kappa/2}/(1+(4-\kappa/2)x_1x_2/(x_2-x_1)^2)$.
Note that, for $\kappa\in(4,8)$:
$$(x_2-x_1)^{4-\kappa/2} \leq r(x_1,x_2)\leq \frac 8\kappa(x_2-x_1)^{4-\kappa/2}.$$
Suppose that the process $(X_1,X_2)$ starts from 0, and define
recursively for $h>0$, $n\geq 0$:
\begin{align*}
S^h_n&=\inf(t\geq T^h_{n-1}: X_{1,t}=X_{2,t}=0)\\
T^h_n&=\inf(t\geq S^h_n: (X_{1,t},X_{2,t})\in M_h)
\end{align*}
with the convention $T^h_{-1}=0$. The number of upcrossings from level
0 to level $h$ at time $t$ is:
$$U(h,t)=\sup(n\in\N:T^h_{n-1}\leq t).$$

\begin{Lem}
The a.s. limit $\lim_n 2^{-n}U(2^{-n},t)$ exists and defines a local
time $(\ell_t)$ for the Markov process $(X_1,X_2)$ at $(0,0)$. The
right-continuous inverse of $\ell$ is a stable subordinator with index
$(2-\kappa/4)$.
\end{Lem}
\begin{proof}
We transpose a classical argument (see \cite{IMK}). Let $(\mu_n)$ be a
decreasing series converging to zero. For a fixed $t>0$, consider
$H_{-n}=\mu_nU(\mu_n,t)$ for $n\in \N$. Then $(H_n)_{n\leq 0}$ is a
${\mc H}$-reversed martingale, where ${\mc H}_n=\sigma(H_m,m\leq n)$;
hence it converges a.s. and in $L^1$ (see e.g. \cite{RW}). Indeed, if
$m<n$, any $\mu_m$-upcrossing contains a $\mu_n$-upcrossing, and the
probability that a given $\mu_n$-upcrossing is contained in a
$\mu_m$-upcrossing is $\mu_n/\mu_m$ (which is the probability that
$R$, started from level $\mu_n$, reaches level $\mu_m$ before
returning to zero). Setting $(\mu_m)=(2^{-m})$, this defines a local
time for $(X_1,X_2)$ at $(0,0)$; in particular $(0,0)$ is a regular
point ($4<\kappa<8$).

The stability index is a straightforward consequence of Brownian
scaling for $(X_1,X_2)$ and the homogeneity property of $r$. More
precisely, let $\lambda>0$; define $\tilde
X_{i,t}=\lambda^{-1}X_{i,\lambda^2t}$, so that $(\tilde X_1,\tilde
X_2)$ and $(X_1,X_2)$ have the same law. Then, if $\tilde U$ designates
the number of upcrossings for $\tilde R=r(\tilde X_1,\tilde X_2)$, one
has
$$\tilde U(h,t)=U(h\lambda^{4-\kappa/2},\lambda^2t).$$
Let $\mu$ be the decreasing sequence such that $\{\mu_m\}_{m\geq
  0}=\{2^{-m}\}_{m\geq 0}\cup\{\lambda^{\kappa/2-4}2^{-m}\}_{m\geq
  0}$. Applying the above result, one gets the a.s. identity:
$\tilde\ell_t=\lambda^{\kappa/2-4}\ell_{\lambda^2 t}.$
Since $\ell$, $\tilde\ell$ are identical in law, this implies that
right-inverse of $\ell$, which is a subordinator by construction, is a
stable subordinator with index $(2-\kappa/4)$.
\end{proof}

By construction, the set of cut-times of the original
$\SLE(\kappa,\kappa-4,\kappa-4)$ (i.e. times at which the trace lies
both on the left and right boundaries of the final hull $K_\infty$) is
the zero set of the process $(X_1,X_2)$. 

\begin{Cor}
The set of cut-times of $\SLE(\kappa,\kappa-4,\kappa-4)$ has
a.s. Hausdorff dimension $(2-\kappa/4)$, $4<\kappa<8$. More precisely, the
$\tilde\varphi_\kappa$-measure of cut-times is a.s. positive and
locally finite, where:
$$\tilde\varphi_\kappa(\eps)=\eps^{2-\kappa/4}(\log|\log\eps|)^{\kappa/4-1}.$$
Moreover, the $\tilde\varphi_\kappa$-measure of cut-times up to time
$t$ is a multiple of the local time of $(X_1,X_2)$ at $(0,0)$.
\end{Cor}
\begin{proof}
Since the right-continuous inverse of $\ell$ is a stable subordinator
with index $(2-\kappa/4)$, the zero set of $(X_1,X_2)$ has the same
distribution as the zero set of a semi-stable process with index
$(\kappa/4-1)^{-1}$, and the result follows from \cite{TW}. 
\end{proof}

Let $\tau$ be the right-continuous of $\ell$. As $\tau$ is stable, one
can formulate excursion-type results w.r.t. cut-times. Note that the
final hull of an $\SLE(6,2,2)$ is identical in law to the filling of a
Brownian excursion (see \cite{D4}, Section 6), since these are two
realizations of the (unique) restriction measure with index 1 (see \cite{LSW3}). So, for
$\kappa=6$, the corollary we presently state connects with results in \cite{Vir}.

\begin{Cor}
(i) The process $(H_u,w_u)_u=(K_{\tau_u},W_{\tau_u})_u$ is a ${\mc
  Q}\times\R$-valued stable L\'evy process, with index
$(2-\kappa/4)$. Moreover, for $\kappa\geq 6$, the following
restriction formula holds: if $A$ is a smooth $+$-hull, $\alpha=\alpha(\kappa,2\kappa-8)=(\kappa-2)(\kappa-3)/2\kappa$, and
$L$ is an independent loop-soup with intensity $\lambda_\kappa$, then : 
$$\phi_A'(0)^\alpha=\E(\phi'_{\phi_{H_u}(A)}(w_u)^\alpha\ind_{H_u\cap A^L=\varnothing}).$$
(ii) Let $e_u=(g_{\tau_{u^-}}(K_{\tau_u}\setminus
K_{\tau_{u^-}}),W_{\tau_u}-W_{\tau_{u^-}})$ if $\tau_u>\tau_{u^-}$,
and $e_u=\partial$ if $\tau_u=\tau_{u^-}$. Then $(e_u)_u$ is a $({\mc
  Q}\times\R)\cup\{\partial\}$-valued, $({\mc X}\vee{\mc Y})_{\tau_u}$-
Poisson process.
\end{Cor}
\begin{proof}
\noindent (i) As in the one-sided case, this follows from the fact that if $T$
is a ${\mc X}\vee{\mc Y}$- stopping time, with $X_{1,T}=X_{2,T}=0$
a.s., then: 
$$(g_{T+s}\circ g_T^{-1}(W_T^-)-W_{T+s},g_{T+s}\circ g_T^{-1}(W_T^+)-W_{T+s},X_{1,T+s},X_{2,T+s})_s$$ 
is distributed as $(Y_{1,t},Y_{2,t},X_{1,t},X_{2,t})_t$ started from
$(0^-,0^+,0^-,0^+)$ and is independent from $({\mc
  X}\vee{\mc Y})_T$ (see proof or Proposition \ref{BP3}). The
restriction formula is then a consequence of \cite{D4}, Section 6.
Assertion (ii) follows
immediately.
\end{proof}

As in the one-sided case, one can conjecture that the L\'evy process
$(K_{\tau_u},W_{\tau_u})_u$ is characterized (up to a scale factor) by
the stability property and the restriction formula (if $\kappa\geq
6$). In the case $\kappa=6$, since $\lambda_6=0$, the law of the total
hull is characterized by the restriction formula.

Finally, we list some remaining questions. It is very likely that
${\mc X}^i\vee{\mc Y}={\mc X}^i$, $i=1,2$, though proving this seems a
bit messy. Getting a constructive (Schramm-Loewner) description of the excursion measure
(``bead measure'') does not look quite straightforward, and may help
to understand the law of $(\tau_u,W_{\tau_u})_u$, where the first
marginal is a stable subordinator with index $(2-\kappa/4)$ and the
second marginal is a symmetric stable process with index $(4-\kappa/2)$. 
A last problem, essentially equivalent to describing the bead measure,
concerns beads ``conditioned to have infinite lifetime''; the first
step consists in taking the Doob $h$-transform of the process $(X_1,X_2)$, using
the ``scale function'' $r$. In
particular, for $\kappa=6$, it should give (for an appropriate
conditioning procedure) a restriction measure with index 2 (see \cite{W2}).

\section{Proof of Watts' formula}

In this section, we present a proof of Watts' formula, that describes
the probability that there exists a double crossing in a rectangle
(top-bottom and left-right) in the scaling limit of critical
percolation. This formula was derived by Watts using (non-rigorous)
Conformal Field Theory techniques (\cite{Wa}). In \cite{D2}, Section 5, we
discussed how Watts' formula could be rephrased in $\SLE_6$
terms. We now sum up this discussion, for the reader's convenience.

Consider a Jordan domain $(D,a,b,c,d)$ with four points marked on the
boundary (in counterclockwise order, say). Suppose that a portion of
the triangular lattice with mesh $\eps$ approximates this domain; each
site of the lattice is colored in blue or yellow with probability
$1/2$, and all sites are independent (this is critical site percolation on the
triangular lattice). Denote by $C_b(A,B)$
(resp. $C_y(A,B)$) the fact that two
site subsets $A$ and $B$ are connected by a blue (resp. yellow)
path, and $T_b(A,B,C)$
(resp. $T_y(A,B,C)$) the fact that $A$, $B$ and $C$ are all connected
by a blue (resp. yellow) cluster of sites. Cardy's formula gives the
probability of events of type
$C_b((ab),(cd))$ in the scaling limit
($\eps\searrow 0)$, and similarly Watts' formula gives the probability
of the event $\{C_b((ab),(cd)),C_b((bc),(da))\}$. For plane topology
reasons, it appears that:
$$
\{C_b((ab),(cd)),C_b((bc),(da))\}=\{T_b((ab),(bc),(cd))\}\setminus\{C_y((ab),(cd))\}$$
$$\{T_b((ab),(bc),(cd))\}=\{C_b((ab),(cd))\}\setminus\{C_b((ab),(cd)),T_y((ab),(bc),(cd))\}.
$$
Since switching the colors of all the sites is a measure-preserving
operation, it follows that:
$$\P(C_b((ab),(cd)),C_b((bc),(da)))=\P(C_b((ab),(cd)))-2\P(T_b((ab),(bc),(cd)),C_y((ab),(cd))).$$

Consider now a Jordan domain $(T,a,b,c)$ with three points marked on
the boundary, and suppose that the sites on $(ab)$ are set to blue and
the sites on $(ca)$ are set to yellow (the arc $(bc)$ remaining
free). Then one can define an exploration process starting from $a$
and stopped when it reaches $(bc)$; this exploration process is the
interface between blue sites connected to $(ab)$ and yellow sites
connected to $(ca)$. As the mesh goes to zero, Smirnov's results
(\cite{Sm1}) imply that this exploration
process converges to chordal $\SLE(6)$ in $(D,a,x)$ stopped when it
hits $(bc)$, where $x$ can be chosen arbitrarily on $(bc)$ (which reflects the
locality property of $\SLE(6)$\,).

Let $X$ be the random endpoint of this exploration process (i.e. $X$
is the first point on $(bc)$ reached by this process). If $(T,a,b,c)$ is
an equilateral triangle, then the distribution of $X$ is uniform on $(bc)$;
this is Carleson's approach  of Cardy's formula. Let $D$ (resp. $E$)
be the lowest point on $(ab)$ (resp. on $(ca)$) reached by the process
before $X$. The exploration hull defines a (random) conformal
quadrilateral $(K,a,D,X,E)$. From the self-duality of the triangular
lattice, this quadrilateral is either crossed by a yellow path
connecting $(aD)$ and $(XE)$ (in which case the exploration process
visits $E$ before $D$) or by a blue path connecting $(DX)$ and $(Ea)$
(and $D$ is visited before $E$). 

Now, let $x$ be some point on $(bc)$. Then, the event $C_b((ab),(xc))$
(with free boundary conditions) is equivalent to $X\in
(xc)$. Moreover, the event $\{T_y((ab),(xc),(ca)),C_b((ab),(xc))\}$ is equivalent to
$X\in (xc)$ and $E$ is visited before $D$. Hence:
\begin{align*}
\P(C_b((ab),(xc)),C_b((bx),(ca)))&=&\P(C_b((ab),(xc)))-2\P(T_y((ab),(xc),(ca)),C_b((ab),(xc)))\\
&=&\P(X\in(xc))-2\P(X\in(xc),E{\textrm{\ visited\ before\ }} D).
\end{align*}

Recall that chordal $\SLE_6$ is the scaling limit of
percolation interfaces for critical site percolation on the triangular
lattice (see \cite{Sm1}). More precisely, consider a chordal $\SLE_6$ in
$(\H,0,\infty)$; here $(T,a,b,c)=(\H,0,1,\infty)$. Then the distribution of $\gamma_{\tau_1}=\inf(\gamma\cap (1,\infty))$
is given by Cardy's formula:
$$\P(1<\gamma_{\tau_1}<x)=B(1/3,1/3)^{-1}\int_{1/x}^1\frac
{ds}{(s(1-s))^{2/3}}=B(1/3,1/3)^{-1}\int_1^x\frac{ds}{(s(s-1))^{2/3}}$$
for any $x>1$. Now, let $g=\sup(t<\tau_1: \gamma_t\in\R)$ be the last
time before $\tau_1$ spent by the trace on the real line. Then Watts'
formula can be rephrased as follows:
$$\P(\gamma_g<0|\gamma_{\tau_1}\in dx)=B(2/3,2/3)^{-1}\int_{1/x}^1\frac{ds}{(s(1-s))^{1/3}}dx.$$
This conditional probability can be derived from the study of
$\SLE(6,2,2)$, as we presently explain. Notations ($d$, $k$,
$\Upsilon$, ...) are as in Section 3.

\begin{Lem} Consider an $\SLE(6,2,2)$ process in $(\H,0,\infty)$, started from
  $(0,y_1,y_2)$, where $y_1<0<y_2$. Let $X_1$ be the leftmost
  swallowed point on $(y_1,0)$ and $X_2$ the rightmost swallowed point
  on $(0,y_2)$. Then the probability that the $\SLE$ reaches $X_1$
  before $X_2$ is given by
$$B(2/3,2/3)^{-1}\int_0^t\frac{ds}{(s(1-s))^{1/3}}$$
where $t=y_2/(y_2-y_1)$.
\end{Lem}
\begin{proof}
The first part of the proof holds for a general value of
$\kappa>4$. We have seen that the distribution of $(X_1,X_2)$ is
$\Upsilon((y_1,y_2),.)$. Let $\tau=\tau_1\wedge\tau_2$, where $\tau_i$ is the first time the trace
reaches $X_i$, $i=1,2$. 
Conditionally on $(X_1,X_2)$, the process
$(g_t(X_1)-W_t,g_t(X_2)-W_t)$, stopped at time $\tau$, is Markov with
semigroup $Q^{\downarrow\downarrow}$. From scale invariance, one can
write 
$$\P_{(x_1,x_2)}^{\downarrow\downarrow}(\tau_1<\tau_2)=f(x_2/(x_2-x_1))$$
for some function $f$ with $f(0)=0$, $f(1)=1$. It is then standard
that $(f(X_{2,t}/(X_{2,t}-X_{1,t})))_{t\geq 0}$ is a martingale (in
the ${\mc X}\vee{\mc Y}$ filtration). From the definition of
$Q^{\downarrow\downarrow}$, it appears that $f$ is such that
$$g'_t(X_1)g'_t(X_2)f(X_{2,t}/(X_{2,t}-X_{1,t}))\partial_{12} k(X_{1,t},X_{2,t})$$
is a $\P$-local martingale, conditionally on $(X_1,X_2)$. So
$(x_1,x_2)\mapsto f(x_2/(x_2-x_1))$ annihilates the differential
operator:
\begin{align*}
\lefteqn{(\partial_{12}k)^{-1}\left(\frac\kappa 2(\partial_1+\partial_2)^2+\frac
  2{x_1}\partial_1+\frac
  2{x_2}\partial_2-\frac 2{x_1^2}-\frac
  2{x_2^2}\right)(\partial_{12}k)=}\\
&\frac\kappa 2(\partial_1+\partial_2)^2+\frac
  2{x_1}\partial_1+\frac
  2{x_2}\partial_2-\left(\frac 1{x_1}+\frac 1{x_2}\right)\frac{4+(4-\kappa)(4-\kappa/2)x_1x_2/(x_2-x_1)^2}{1+(4-\kappa/2)x_1x_2/(x_2-x_1)^2}(\partial_1+\partial_2).
\end{align*}
It follows that $f$ solves the ODE:
$$f''(t)+\left(-\frac 4\kappa\left(\frac 1t+\frac 1{t-1}\right)+2\frac
  {(4-\kappa/2)(2t-1)}{1+(4-\kappa/2)t(t-1)}\right)f'(t)=0$$
which implies that
$$f(t)=c\int_0^t\frac{(s(1-s))^{4/\kappa}}{(1+(4-\kappa/2)s(s-1))^2}ds$$
where $c$ is such that $f(1)=1$. Let 
$$h(y_1,y_2)=\P_{(y_1,y_2)}^{\uparrow\uparrow}(\tau_1<\tau_2).$$
Since the distribution of $(X_1,X_2)$ is $\Upsilon((y_1,y_2),.)$, and
the conditional probability of $\{\tau_1<\tau_2\}$ given $(X_1,X_2)$ is
$f(X_2/(X_2-X_1))$, one gets, integrating by parts:
\begin{align*}
h(y_1,y_2)&=\int_{y_1}^0\int_0^{y_2}f\left(\frac{x_2}{x_2-x_1}\right)\frac{-\partial_{12}
  k(x_1,x_2)}{k(y_1,y_2)}dx_1dx_2\\
&=\int_{y_1}^0 \left(\frac{-\partial_1
    k(x_1,y_2)}{k(y_1,y_2)}f\left(\frac{y_2}{y_2-x_1}\right)-\int_0^{y_2}\frac{-\partial_1 k(x_1,x_2)}{k(y_1,y_2)}\frac{-x_1}{(x_2-x_1)^2}f'\left(\frac{x_2}{x_2-x_1}\right)dx_2\right)dx_1\\
&=f\left(\frac
  {y_2}{y_2-y_1}\right)+\int_{y_1}^0\frac{k(x_1,y_2)}{k(y_1,y_2)}f'\left(\frac{y_2}{y_2-x_1}\right)\frac{y_2}{(y_2-x_1)^2}dx_1\\
&\hphantom{==}-\int_{y_1}^0\int_0^{y_2}\frac{\partial_1 k(x_1,x_2)}{k(y_1,y_2)}\frac{x_1}{(x_2-x_1)^2}f'\left(\frac{x_2}{x_2-x_1}\right)dx_1dx_2\\
&=f\left(\frac
  {y_2}{y_2-y_1}\right)
+c\int_{y_1}^0\frac{(-x_1)y_2^{1+4/\kappa}(y_2-x_1)^{\kappa/2-6}}{y_1^d(y_2-y_1)^{\kappa
    d^2/2}(1+(4-\kappa/2)(x_1y_2)/(y_2-x_1)^2)^2}dx_1\\
&\hphantom{==}+cd\int_{y_1}^0\int_0^{y_2}\left(\frac
  1{x_1}+\frac\kappa 2\frac d{x_1-x_2}\right)x_1^2x_2\frac{(x_2-x_1)^{\kappa/2-6}}{(1+(4-\kappa/2)x_1x_2/(x_2-x_1)^2)^2}\frac{dx_1dx_2}{k(y_1,y_2)}. 
\end{align*}
Now, let $\kappa=6$. The two integrals in the former expression have
rational integrands, so they can be computed mechanically. After
simplifications, one gets:
$$h(y_1,y_2)=f\left(\frac
  {y_2}{y_2-y_1}\right)+\frac c3\frac
{y_1y_2(y_1+y_2)}{(y_1^2-y_1y_2+y_2^2)(-y_1y_2)^{1/3}(y_2-y_1)^{1/3}}.$$
By homogeneity, $h(y_1,y_2)=\tilde h(t)$, where
$t=y_2/(y_2-y_1)$. Differentiating the above expression:
\begin{align*}
c^{-1}\tilde h'(t)&=\frac{(t(1-t))^{2/3}}{(1+t(t-1))^2}+\frac
d{dt}\left(\frac {(1-2t)(t(1-t))^{2/3}}{3(1+t(t-1))}\right)\\
&=\frac 29 (t(1-t))^{-1/3}
\end{align*}
Hence $c^{-1}=2B(2/3,2/3)/9$, and 
$$\P_{(y_1,y_2)}^{\uparrow\uparrow}(\tau_1<\tau_2)=B(2/3,2/3)^{-1}\int_0^t\frac{ds}{(s(1-s))^{1/3}}$$
where $t=y_2/(y_2-y_1)$.
\end{proof}

Consider now critical site percolation on the triangular lattice. In
this case, interfaces converge to $\SLE_6$ in the
scaling limit (see \cite{Sm1}).

\begin{Prop}[Watts' formula]
The probability that the four boundary arcs of a conformal
quadrilateral $(D,a,b,c,d)$, where $D$ is a simply connected Jordan
domain, are connected by a cluster is:
$$B(1/3,1/3)^{-1}\int_0^z\frac{1}{(s(1-s))^{2/3}}\left(1-2B(2/3,2/3)^{-1}\int_0^s\frac{dr}{(r(1-r))^{1/3}}\right)ds$$
where $z$ is the cross-ratio $[a,b,c,d]$.    
\end{Prop}
\begin{proof}
Consider a chordal $\SLE(6)$ in $(\H,0,\infty)$; $\gamma$ is the
trace, $\tau_1$ is the swallowing time of $1$, and $g$ is the last
time before $\tau_1$ spent by the trace on $\R$. As we have seen,
$(1/\gamma_{\tau_1})$ has distribution $\Beta(1/3,1/3)$ (Cardy's
formula), and Watts' formula is equivalent to
$$\P(\gamma_g<0|\gamma_{\tau_1}\in
dx)=B(2/3,2/3)^{-1}\int_{1/x}^1\frac{ds}{(s(1-s))^{1/3}}dx.$$
So let $(x_1,x_2)$ be a neighbourhood of $x>1$. From the locality
property of $\SLE(6)$ (see e.g. \cite{W1}), the $\SLE(6)$
conditionally on $\gamma_{\tau_1}\in (x_1,x_2)$ stopped at
$\tau_{x_1}$ is identical in law to a time-changed chordal $\SLE(6)$ in $(\H,0,x)$
with the corresponding conditioning, stopped when its trace hits
$(x_1,x_2)$. Consider the homography $\phi(z)=z(1-x)/(z-x)$. Then the
image of this conditioned $\SLE$ under $\phi$ is a chordal $\SLE(6)$
in $(\H,0,\infty)$, conditioned on its trace not hitting
$(y_2,1-x)\cup(1,y_1)$, where $y_i=\phi(x_i)$, $i=1,2$. What we have
to compute is the probability that the trace is on $(1-x,0)$ the last
time it visits $(1-x,1)$.

As $x_1\nearrow x$, $x_2\searrow x$, it appears that
$y_1\rightarrow\infty$, $y_2\rightarrow -\infty$; formally, we have now
an $\SLE(6)$ conditioned on its trace not hitting $(-\infty,1-x)\cup
(1,\infty)$. As we have seen, this singular conditioning can be
realized as an $\SLE(6,2,2)$ process started from
$(0,1-x,1)$. According to the previous lemma, the probability that the
trace is on $(1-x,0)$ the last time it visits $(1-x,1)$, for an
$\SLE(6,2,2)$ process, is:
$$B(2/3,2/3)^{-1}\int_{1/x}^1\frac{ds}{(s(1-s))^{1/3}}$$
which is exactly what we need. To make the previous limiting argument
more precise, one can argue along the lines of Theorem 3.1 in
\cite{LSW6}, as we now sketch. If $(g_t)$ denotes the family of
conformal equivalences associated with a chordal $\SLE(\kappa)$
process, $\kappa>4$, and $W$ is its driving process, let
$$Z_t=\frac{W_t-g_t(u_1)}{g_t(u_2)-g_t(u_1)} $$
for some $u_1<0<u_2$. Then, after an appropriate time-change, $Z$ is a diffusion on
$[0,1]$ with leading eigenvector $h(z)=(z(1-z))^{1-4/\kappa}$ (see \cite{LSW6}). Using
this eigenvector, one can ``condition'' this diffusion to have infinite
lifetime (i.e. $Z$ never swallows $0$ or $1$, that is the $\SLE$ never
swallows $u_1$ or $u_2$). Working backwards, it appears that the
corresponding conditional $\SLE$ is precisely
$\SLE(\kappa,\kappa-4,\kappa-4)$ started from $(0,u_1,u_2)$. Suppose
now that $M\gg 1$; if a chordal $\SLE(\kappa)$ run until
time $M^2$ has not swallowed $u_1$, $u_2$ yet, then the conditional
probability that the trace does not hit $(-M,u_1)\cup(u_2,M)$ is
bounded away from $0$. Conversely, if the trace does not hit
$(-M,u_1)\cup(u_2,M)$, $M\rightarrow\infty$, then the lifetime of $Z$
goes to infinity. So the two procedures - conditioning on the
trace not hitting $(-M,u_1)\cup (u_2,M)$, $M\rightarrow\infty$, or
conditioning on $Z$ not hitting $0$ or $1$ before time $T$,
$T\rightarrow\infty$\ - yield the same limiting object, namely $\SLE(\kappa,\kappa-4,\kappa-4)$. 

Alternatively, one can build on the two following facts (for
$4<\kappa<8$):
\begin{itemize}
\item Let $0<x<y$, $\tau_x$ the swallowing time of $x$. A chordal
  $\SLE(\kappa)$ process in $(\H,0,\infty)$ conditioned
  on $\gamma_{\tau_x}\in dy$ and stopped at $\tau_y$ is an $\SLE(\kappa,\kappa-4,-4)$ process
  in $(\H,0,\infty)$ started from $(0,x,y)$ and stopped at $\tau_y$.
\item A chordal $\SLE(\kappa,\rho_1,\rho_2)$ in $(\H,0,\infty)$
  started from $(0,x,y)$ is a time-changed
  $\SLE(\kappa,\kappa-6-\rho_1-\rho_2,\rho_1)$ in $(\H,0,y)$ started
  from $(0,\infty,x)$.
\end{itemize}
The first fact is a direct consequence of Cardy's formula for $\SLE$
and Girsanov's theorem. The second fact can be proved by a computation along
the lines of \cite{LSW3}, Section 5. Consequently, an $\SLE(\kappa)$ conditioned on $\gamma_{\tau_x}\in dy$
is a time-changed $\SLE(\kappa,2,\kappa-4)$ in $(\H,0,y)$ started from
$(0,\infty,y)$. For $\kappa=6$, one recovers the required result.
\end{proof}

Originally, Watts identified this function as the only function
annihilating the
fifth-order differential operator
$$(z(1-z))^{-2}\frac{d^3}{dz^3}(z(1-z))^{4/3}\frac{d}{dz}(z(1-z))^{2/3}\frac{d}{dz}$$
and satisfying appropriate boundary conditions. This operator seems to
have no obvious interpretation in the $\SLE_6$ framework. As pointed
out in \cite{KZ}, the solution space of the associated ODE is spanned
by: the constant function $\ind$, the function appearing in Cardy's
formula, the function appearing in Watts' formula, $z\mapsto\log(z)$,
$z\mapsto\log(1-z)$. Moreover, Cardy's prediction for the expected
number of disjoint clusters connecting two opposite sides of a
rectangle also belongs to this solution space (\cite{Ca1}, see also
\cite{Maier}). 

It is also possible to express Watts' formula using equianharmonic
elliptic functions. From the known value of the equianharmonic
$\sigma$ function at half-periods, one can deduce that Watts' formula
imply the following result (\cite{Maier}).

\begin{Cor}
The probability that there exists a double crossing in a square is
$$\frac 14+\frac{\sqrt 3}{4\pi}(3\log 3-4\log 2)\simeq 0.322120455\dots$$
\end{Cor}

Note that the numerical value of this probability agrees with the estimate in \cite{LPS}.

\bibliographystyle{abbrv}
\bibliography{biblio}

\begin{thebibliography}{10}

\bibitem{Bef}
V.~Beffara.
\newblock {Hausdorff dimensions for $\SLE_6$}.
\newblock {\em preprint, arXiv:math.PR/0204208}, 2002.

\bibitem{BerComp}
J.~Bertoin.
\newblock Complements on the {H}ilbert transform and the fractional derivative
  of {B}rownian local times.
\newblock {\em J. Math. Kyoto Univ.}, 30(4):651--670, 1990.

\bibitem{Ber90}
J.~Bertoin.
\newblock Excursions of a {${\rm BES}\sb 0(d)$} and its drift term {$(0<d<1)$}.
\newblock {\em Probab. Theory Related Fields}, 84(2):231--250, 1990.

\bibitem{Ca1}
J.~L. Cardy.
\newblock Conformal invariance and percolation.
\newblock {\em preprint, arXiv:math-ph/0103018}, 2001.

\bibitem{CPY}
P.~Carmona, F.~Petit, and M.~Yor.
\newblock Beta-gamma random variables and intertwining relations between
  certain {M}arkov processes.
\newblock {\em Rev. Mat. Iberoamericana}, 14(2):311--367, 1998.

\bibitem{D2}
J.~Dub{\'e}dat.
\newblock Reflected brownian motions, intertwining relations and crossing
  probabilities.
\newblock {\em Ann. Inst. H. Poincar\'e Probab. Statist., to appear}, 2003.

\bibitem{D4}
J.~Dub{\'e}dat.
\newblock $\sle(\kappa,\rho)$ martingales and duality.
\newblock {\em Ann. Probab., to appear}, 2003.

\bibitem{IMK}
K.~It{\^o} and H.~P. McKean, Jr.
\newblock {\em Diffusion processes and their sample paths}.
\newblock Springer-Verlag, Berlin, 1974.
\newblock Second printing, corrected, Die Grundlehren der mathematischen
  Wissenschaften, Band 125.

\bibitem{KZ}
P.~Kleban and D.~Zagier.
\newblock Crossing probabilities and modular forms.
\newblock {\em J. Statist. Phys.}, 113(3-4):431--454, 2003.

\bibitem{LPS}
R.~Langlands, P.~Pouliot, and Y.~Saint-Aubin.
\newblock Conformal invariance in two-dimensional percolation.
\newblock {\em Bull. Amer. Math. Soc. (N.S.)}, 30(1):1--61, 1994.

\bibitem{LSW3}
G.~Lawler, O.~Schramm, and W.~Werner.
\newblock Conformal restriction: the chordal case.
\newblock {\em J. Amer. Math. Soc.}, 16(4):917--955 (electronic), 2003.

\bibitem{LSW2}
G.~F. Lawler, O.~Schramm, and W.~Werner.
\newblock Conformal invariance of planar loop-erased random walks and uniform
  spanning trees.
\newblock {\em Ann. Inst. H. Poincar\'e Statist. Probab., to appear.}, 2002.

\bibitem{LSW6}
G.~F. Lawler, O.~Schramm, and W.~Werner.
\newblock Values of {B}rownian intersection exponents. {III}. {T}wo-sided
  exponents.
\newblock {\em Ann. Inst. H. Poincar\'e Probab. Statist.}, 38(1):109--123,
  2002.

\bibitem{Maier}
R.~S. Maier.
\newblock On crossing event formulas in critical two-dimensional percolation.
\newblock {\em J. Statist. Phys.}, 111(5-6):1027--1048, 2003.

\bibitem{PY84}
J.~Pitman and M.~Yor.
\newblock A decomposition of {B}essel bridges.
\newblock {\em Z. Wahrsch. Verw. Gebiete}, 59(4):425--457, 1982.

\bibitem{RY}
D.~Revuz and M.~Yor.
\newblock {\em Continuous martingales and {B}rownian motion}, volume 293 of
  {\em Grundlehren der Mathematischen Wissenschaften}.
\newblock Springer-Verlag, Berlin, third edition, 1999.

\bibitem{RW}
L.~C.~G. Rogers and D.~Williams.
\newblock {\em Diffusions, {M}arkov processes, and martingales. {V}ol. 1}.
\newblock Wiley Series in Probability and Mathematical Statistics: Probability
  and Mathematical Statistics. John Wiley \& Sons Ltd., Chichester, second
  edition, 1994.
\newblock Foundations.

\bibitem{RS01}
S.~Rohde and O.~Schramm.
\newblock {Basic Properties of SLE}.
\newblock {\em Ann. Math., to appear}, 2001.

\bibitem{S0}
O.~Schramm.
\newblock Scaling limits of loop-erased random walks and uniform spanning
  trees.
\newblock {\em Israel J. Math.}, 118:221--288, 2000.

\bibitem{Sm1}
S.~Smirnov.
\newblock {Critical percolation in the plane. I. Conformal Invariance and
  Cardy's formula II. Continuum scaling limit}.
\newblock {\em in preparation}, 2001.

\bibitem{TW}
S.~J. Taylor and J.~G. Wendel.
\newblock The exact {H}ausdorff measure of the zero set of a stable process.
\newblock {\em Z. Wahrscheinlichkeitstheorie und Verw. Gebiete}, 6:170--180,
  1966.

\bibitem{Vir}
B.~Vir{\'a}g.
\newblock Brownian beads.
\newblock {\em Probab. Theory Related Fields}, 127(3):367--387, 2003.

\bibitem{Wa}
G.~M.~T. Watts.
\newblock A crossing probability for critical percolation in two dimensions.
\newblock {\em J. Phys. A}, 29(14):L363--L368, 1996.

\bibitem{W1}
W.~Werner.
\newblock {Lectures on random planar curves and Schramm-Loewner evolution}.
\newblock In {\em Lecture Notes of the 2002 St-Flour summer school, to appear}.
  Springer-Verlag, 2002.

\bibitem{W2}
W.~Werner.
\newblock Girsanov transformation for $\sle(\kappa,\rho)$ processes,
  intersection exponents and hiding exponent.
\newblock {\em Ann. Fac. Sci. Toulouse}, 2003.

\end{thebibliography}

-----------------------

Laboratoire de Math\'ematiques, B\^at. 425

Universit\'e Paris-Sud, F-91405 Orsay cedex, France

julien.dubedat@math.u-psud.fr

\end{document}